\documentclass[runningheads]{llncs}

\usepackage{graphicx}

\usepackage{varwidth}

\usepackage[T2A]{fontenc}
\usepackage[utf8]{inputenc}
\usepackage[russian,english]{babel}

\usepackage{amsmath}
\usepackage{amssymb}
\usepackage{mathtools}
\usepackage{nicefrac}
\usepackage{algorithm,algorithmic}

\usepackage{subfigure}
\usepackage{float}

\usepackage{booktabs}

\usepackage{tabularx}
\usepackage{multirow}
\usepackage{multicol}
\usepackage{makecell}
\usepackage{enumitem}

\usepackage{colortbl}
\usepackage{xcolor}

\usepackage{hyperref}
\hypersetup{colorlinks=true,allcolors=blue}
\usepackage{hypcap}
\usepackage[autostyle]{csquotes}
\usepackage[framemethod=tikz]{mdframed}

\newcommand{\bad}[1]{\colorbox{red}{\color{white}{#1}}}
\newcommand{\ok}[2]{\colorbox{blue!#1}{\color{white}{#2}}}

\DeclareMathOperator*{\argmin}{arg\,min}

\newcommand{\CC}{\mathbb{C}}

\newcommand{\R}{{\mathbb R}} 
\def\R{\mathbb R}

\def\BI{\begin{itemize}}
\def\EI{\end{itemize}}

\def\ba{\begin{array}}
\def\ea{\end{array}}
\def\beann{\begin{eqnarray*}}
\def\eeann{\end{eqnarray*}}
\def\bea{\begin{eqnarray}}
\def\eea{\end{eqnarray}}

\def\diag{{\rm diag \,}}

\def\beq{\begin{equation}}
\def\eeq{\end{equation}}

\def\vf{\varphi}

\def\epsilon{\varepsilon}

\def\BMP{\begin{minipage}{9.5cm}}
\def\EMP{\end{minipage}}
\def\MPT{\begin{minipage}{12.5cm}}
\def\EPT{\end{minipage}}

\def\la{\langle}
\def\ra{\rangle}

\begin{document}

\title{Non-convex optimization in digital pre-distortion of the signal\thanks{The work was supported by the Russian Science Foundation (project 21-71-30005).
}}

\author{Dmitry Pasechnyuk \inst{1}\orcidID{0000-0002-1208-1659}\and
Alexander Maslovskiy \inst{1}\orcidID{0000-0001-7388-6146}\and
Alexander Gasnikov\inst{1,4,5}\orcidID{0000-0002-7386-039X}\and
Anton Anikin\inst{2}\orcidID{0000-0002-7681-2481}\and
Alexander Rogozin\inst{1}\orcidID{0000-0003-3435-2680}\and
Alexander Gornov\inst{2}\orcidID{0000-0002-8340-5729}\and
Andrey Vorobyev \inst{3}\and
Eugeniy Yanitskiy \inst{3}\and
Lev Antonov \inst{3}\orcidID{0000-0002-2403-9503
}\and
Roman Vlasov \inst{3}\orcidID{0000-0002-3798-8310} \and Anna Nikolaeva\inst{6}\orcidID{0000-0003-2980-2587}
\and Maria Begicheva \inst{6}\orcidID{0000-0002-9824-0373}}

\authorrunning{D.~Pasechnyuk et al.}

\institute{
Moscow Institute of Physics and Technology, Dolgoprudny, Russia\\
\email{\{pasechniuk.da,aleksandr.maslovskiy,aleksandr.rogozin\}@phystech.edu}\and
Matrosov Institute for System Dynamics and Control Theory, Irkutsk, Russia\\
\email{\{anikin,gornov\}@icc.ru}\and
Huawei Russian Research Institute, Moscow, Russia\\
\email{\{andrey.vorobyev,yanitskiy.eugeniy,vlasov.roman,antonov.lev\}@huawei.com}
\and Institute for Information Transmission Problems RAS\and Russia
Caucasus Mathematical Center, Adyghe State University, Russia\\ \email{gasnikov@yandex.ru}
\and Skolkovo Institute of Science and Technology \\ \email{\{anna.nikolaeva,maria.begicheva\}@skoltech.ru}
}

\maketitle
\begin{abstract}

In this paper we give some observation of applying modern optimization methods for functionals describing digital predistortion (DPD) of signals with orthogonal frequency division multiplexing (OFDM) modulation.
The considered family of model functionals is determined by the class of cascade Wiener--Hammerstein models, which can be represented as a computational graph consisting of various nonlinear blocks.
To assess optimization methods with the best convergence depth and rate as a properties of this models family we multilaterally consider modern techniques used in optimizing neural networks and numerous numerical methods used to optimize non-convex multimodal functions.

The research emphasizes the most effective of the considered techniques and describes several useful observations about the model properties and optimization methods behavior. 

\end{abstract}

\section{Introduction}
{
Today, base stations, which perform in the capacity of radio signal transceivers, are widely used for the implementation and organization of wireless communication between remote devices. Modern base stations have a complex technical structure and include many technical components allowing organize accurate and efficient data transmission. One of the most important of these components is the analog power amplifier (PA). Its role is to amplify the signal from the base station, reduce the noise effect on the signal and increase the transmission range.



The impact of some ideal amplifiers can be characterized with functional $\mathcal{PA}(x)=a \cdot x$, where $a \gg 1$ and $x$ is an input signal. 
However, real amplifiers are complex non-linear analog devices that couldn't be described by analytical function due to the influence of many external and internal obstructing factors. Power Amplifiers can change phase, cut amplitude of the original signal, and generate parasitic harmonics outside the carrier frequency range. These influences cause significant distortions of the high-frequency and high-bandwidth signal. The spectrum plot from fig.~\ref{fig:spectrum} shows that described problem is relevant in the conditions of operation of modern devices: when the signal goes through the power amplifier its spectrum range becomes wider than the spectrum range of the original signal and as a result generates noise for other signals.

\vspace{-0.5cm}
\begin{figure}[H]
\centering
\vspace{-0.2cm}
\includegraphics[width=0.60\linewidth]{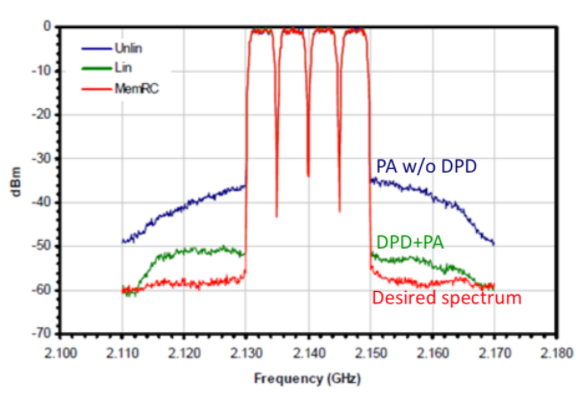}
\vspace{-0.4cm}
\caption{Power spectral density plot of original signal, out of PA signal, and result of pre-distorting signal}
\label{fig:spectrum}
\end{figure}
\vspace{-0.5cm}

One possible solution to this problem is employing the digital baseband pre-distortion (DPD) technique to compensate for non-linear effects, that influence the input signal. In this case, DPD acts upon the input signal with an inverse non-linear with the aim of offsetting the impact of the power amplifier.

\vspace{-0.5cm}
\begin{figure}[H]
\centering
\includegraphics[width=0.60\linewidth]{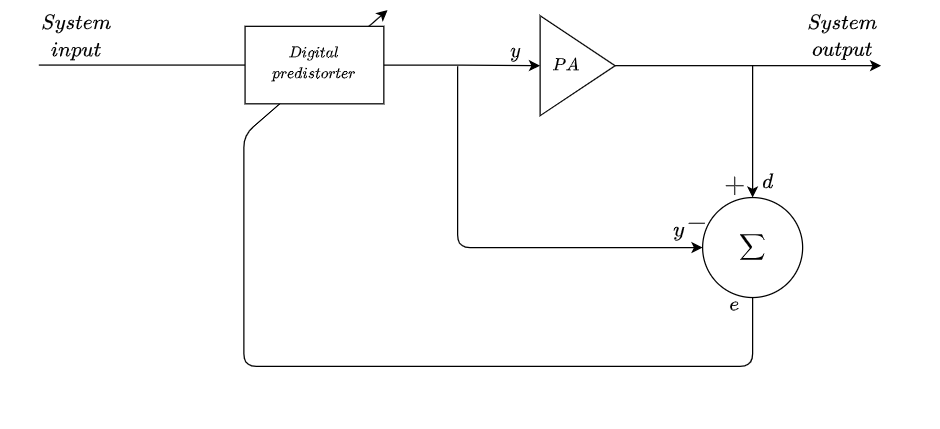}
\vspace{-0.2cm}
\caption{DPD principle of suppression spread spectrum \cite{haykin2008adaptive}}
\label{fig:DPD_process}
\end{figure}
\vspace{-0.5cm}
 
\noindent To compensate for errors generated during the amplification, the difference between input and output of power amplifier signals transmits to adaptive digital pre-distorter to optimize its parameters. In this paradigm, the pre-distorter model can be presented as a parametric function transforming signal in accordance with the real digital pre-distorter operation. Thus, the parametrization of the model is being optimized as parameters of real adaptive filters.





\vspace{-0.5cm}
\begin{figure}[H]
    \centering
    \vspace{-0.2cm}
    \includegraphics[width=0.75 \linewidth]{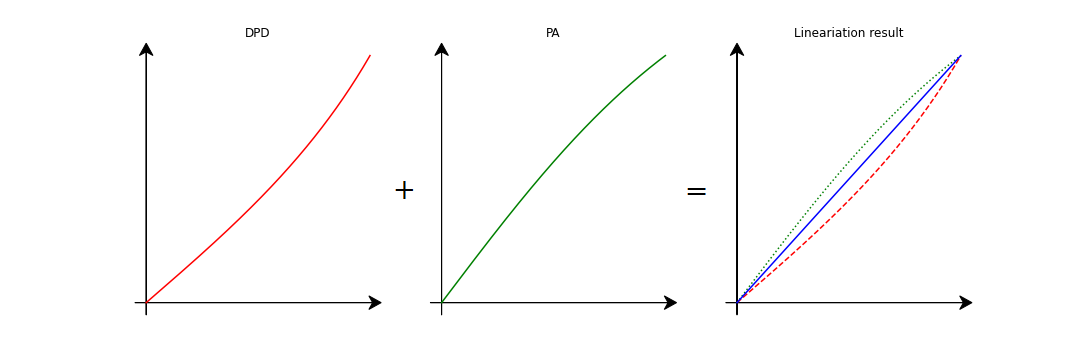}
    \vspace{-0.5cm}
    \caption{DPD and PA influence of signal}
\end{figure}
\vspace{-0.8cm}


According to all the above, from a mathematical point of view, pre-distortion consists of applying to the input signal function $\mathcal{DPD}$, that approximates the inverse to the real $\mathcal{PA}$ function describing the effect of an analog power amplifier \cite{haykin2008adaptive}. This problem can be expressed in the following optimization form\footnote{Note that even in this formulation, the problem can be solved by a classical gradient descent scheme $y^{k+1} = y^k - h (\mathcal{PA}(y^k) - ax)$, assuming that DPD can model arbitrary function, and that the Jacobian $\partial \mathcal{PA}(y) / \partial y \approx a I$.}: 
\[
    \frac{1}{2}\|\mathcal{PA}(y) - ax\|_2^2 \to\min_{y=:\mathcal{DPD}(x)}.
\]

A more practical approach is to choose a certain parametric family of functions $\{DPD(x,\theta)\}_{\theta \in \Theta}$ (in particular, defined by a computational graph of a specific type). Taking into account that non-linear transformation of the signal can be obtained as a result of passing it through a number of non-linear functions and thus presented as the composition of some additive changes, the optimization problem reformulated in the following form:
\[
    \frac{1}{2}\|\mathcal{DPD}_\theta(x) - e\|_2^2 \to\min_{\theta \in \Theta}.
\]

If we have large enough training set\footnote{In fact, it is enough to have just a large segment of the input signal, since by solving one of the problems posed above by a method tuned to a sufficiently high accuracy, after some time it is possible to obtain a sufficiently accurate approximation of the inverse function, which can, in turn, act as a benchmark when carrying out the experiments.} $(x,e)$, it is possible to optimize the model parameters on it, thereby choosing a good approximation for the function that acts as a DPD for the sample signal. 

This work is devoted to a wide range of issues related to the numerical solution of problems of this kind for one fairly wide class of models~--- Wiener--Hammerstein models \cite{ghannouchi2015behavioral,schreurs2008rf}. Based on the results of numerous computational experiments, there were identified and are now described the methods that demonstrate the most successful results in terms of the convergence depth, its rate, and method's susceptibility to overfitting. Approaches to online and offline training of DPD models, methods of initializing models, and also some directions for possible further development of methods for solving the problems of the category under consideration are proposed.

}

\section{Problem formulation}
{

\subsection{Model Description} 
\label{model_subs}

In this paper we consider block oriented models describing dynamic nonlinear effects of PA. These models have a tend to reduce number of coefficients unlike the Volterra series. Considered Wiener--Hammerstein models can account for static non-linear behavior of PA and deal with linear memory effects in the modeled system. To refine model robustness and enhance its performance it was chosen an instance of cascade Wiener--Hammerstein model \cite{ghannouchi2015behavioral}, those structure is presented in fig. \ref{fig:block_model}. 


\begin{figure}[H]
\centering
\includegraphics[width=1\textwidth]{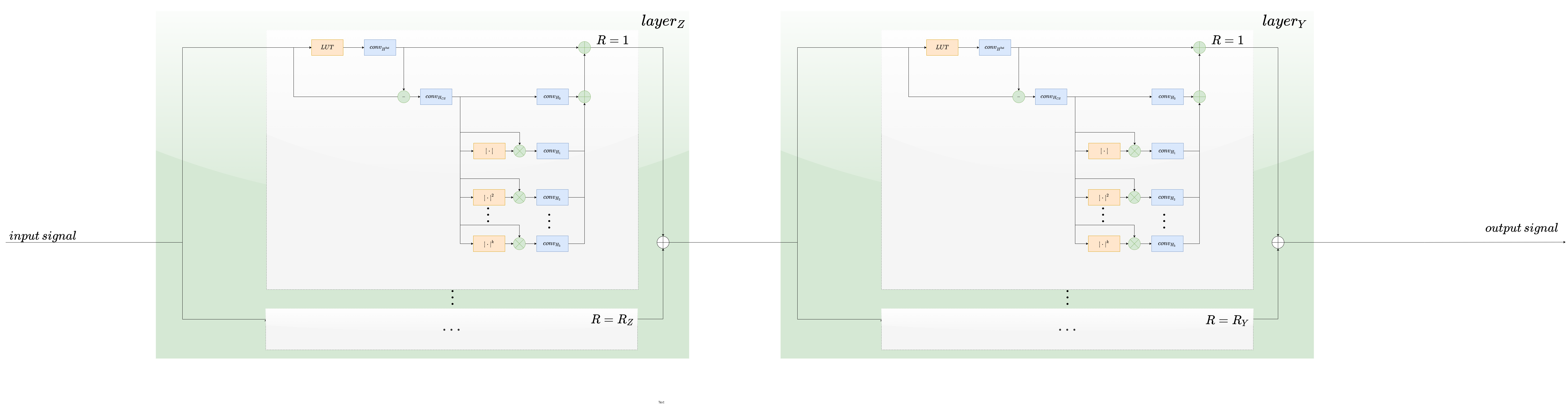}
\caption{Two-layer block model of the Wiener--Hammerstein type}
\label{fig:block_model}
\end{figure}

Let us describe a formal, with some generality, mathematical model for the case of two layers. The first $z$-layer and the second $y$-layer's outputs are the sum of the results returned from $R_z$ and $R_y$, respectively, identical blocks forming the previous layer. Each of these blocks is described as a combination of convolution, polynomials and lookup table functions applied to the input signal:
\[
d_{k,H^{lut,H_{CS}}}(x):=
\text{conv}_{H_{CS}}\left(\text{conv}_{k,H^{lut}}\left(\sum_{p=1}^P C_p \cdot \phi_p(|x|)\cdot x \right)-x\right)
\]
\begin{align}
    \text{block}_{H,H_{CS}, H^{lut}, C, k}(x) &:= \text{conv}_{k,H^{lut}}\left(\sum_{p=1}^P C_p \cdot \phi_p(|x|)\cdot x \right)\\
    \nonumber&+\sum_{l=0}^{B_{k}}\text{conv}_{k,H_{l}}\left(d_{k,H^{lut},H_{CS}}(x)\cdot|d_{k,H^{lut},H_{CS}}(x)|^l\right) 
\end{align}
where $H^{lut} \in \CC^M, H \in \CC^N, H_{CS} \in \CC^{K\times L}=\left\{H_l \in \CC^K\;\big|\;l\in\{1,...,L\}\right\}$ denote weights of convolutions, $C \in \CC^P$ are weighting coefficients of gains in lookup table functions, $\phi_p $ is a polynomial function of arbitrary order applied to the input vector to activate special gain for quantized amplitude of the complex input, $\text{conv}_{k, H} $ is the convolution of the input vector with a vector of weights $H$ and shift $k \in \mathbb {N}$ \cite{ghannouchi2015behavioral}:
\[
    \text{conv}_{k, H}(x) := \sum_{n=1}^{N} H_n x_{k-n+1}. 
\]

Thus, the presented two-layer model is described as follows:
\begin{align}
    z_k(x) &= \sum_{r=1}^{R_z}\text{block}_{H_z,H_{CS,z}, H^{lut}_z, C_z, k}(x), \\
    \nonumber y_k(x) &= \sum_{r=1}^{R_y}\text{block}_{H_y,H_{CS,y}, H^{lut}_y, C_y, k}(z(x)).
\end{align}
In this work, we also study the following modification of the described model, obtained by utilizing skip connections technique (widely used in residual neural networks \cite{he2016deep}):
\begin{align} \label{model}
    z_k(x) &= \sum_{r=1}^{R_z}\text{block}_{H_z,H_{CS,z}, H^{lut}_z, C_z, k}(x), \\
    \nonumber y_k(x) &= \sum_{r=1}^{R_y}\text{block}_{H_y,H_{CS,y}, H^{lut}_y, C_y, k}(z(x)) + z_k(x).
\end{align}

As a result, we get a computational graph, characterized by the following hyperparameters (such a set for each layer): 
\begin{center}
\begin{enumerate}
\item  $N$, $M$, $K$~--- width of applied convolution,
\item  $P$~--- number of spline functions,
\item  $R$~--- number of blocks in a layer;
\end{enumerate}
\end{center}
and having the following set of training parameters:
\[
    \theta := (H_z, H^{lut}_z, H_{CS,z}, H_y, H^{lut}_y,H_{CS,y},C_z, C_y),
\]
where
\begin{align*}
    &\begin{rcases}
    H_z\in \CC^{R_z\times K_z\times L}, H_{CS,z}\in \CC^{R_z\times N_z}, H^{lut}_z\in \CC^{R_z\times M_z} \\
    H_y\in \CC^{R_y\times K_y\times L}, H_{CS,y}\in \CC^{R_y\times N_y}, H^{lut}_y\in \CC^{R_y\times M_y}
    &\end{rcases} \text{  convolution weights} \\
    &\begin{rcases}
    C_z\in \CC^{R_z\times P_z}, C_y\in \CC^{R_y\times P_y}
    &\end{rcases} \text{  spline weights.}
\end{align*}

The total number of model parameters can be calculated as follows: $n = R_z(N_z + M_z + K_z\times L + P_z) + R_y(N_y + M_y + K_y\times L + P_y)$. In the numerical experiments presented in this paper, the tuning of the model was such that the number of model parameters $n \sim 10^3$. Note that additional experiments on considering various graph configurations and hyperparameter settings are presented in Appendix~\ref{tune}.

\subsection{Optimization Problem Statement}

Let's denote the result returned by the used model parameterized by the vector $\theta$ at the input $x$ as $\mathcal{M}_\theta (x) := y(x)$ \eqref{model}. The main considered problem of restoring the function $\mathcal{PA}^{-1}$ using the described model can be formulated as a problem of supervised learning in the form of regression. Let then $(x, \overline{y})$ be a training sample, where $x \in \CC^m$ is a signal input to the DPD, $\overline{y} \in \CC^m$ is the desired modulated output signal. In this setting, the problem of restoring the $\mathcal{DPD}$ function can be formulated as minimizing the empirical risk (in this case, with a quadratic loss function):
\begin{equation} \label{main_problem}
    f(\theta) := \frac{1}{m} \sum_{k=1}^m ([\mathcal{M}_\theta (x)]_k - \overline{y}_k)^2 \rightarrow \min_\theta.
\end{equation}
To assess the quality of the solution obtained as a result of the optimization of this loss functional, we will further use the normalized mean square error quality metric, measured in decibels: 
\[
    \text{NMSE}(y, \overline y) := 10 \log_{10}\left\{\frac{\sum_{k=1}^{m} (y_k - \overline y_k)^2}{\sum_{k=1}^m x_k^2}\right\} \quad \text{dB}.
\]

}

\section{Optimization Methods}
{
In this main section of the article, we consider three wide classes of optimization methods: full-gradient methods, Gauss--Newton methods, and stochastic (SGD-like) methods. Only the most effective methods for the problem under consideration are described directly, however, the provided experimental results cover a much wider variety of algorithms. The implementation features of some of the approaches are also described.

\subsection{Long memory L-BFGS}

Let us start with considering the class of quasi-Newton methods. In contrast to the classical Newton's method, which uses the Hessian to find the quadratic approximation of a function at a certain point, the quasi-Newton methods are based on the principle of finding such a quadratic approximation that is tangent to the graph of the function at the current point and have the same gradient value as the original function at the previous point of the trajectory. More specifically, these methods have iterations of form $x_{k+1} = x_k - h_k H_k 
\nabla f(x_k)$, where $H_k$ is an approximation of inverse Hessian $[\nabla^2 f(x_k)]^{-1}$ and $h_k$ is the step-size. The choice of matrices $H_k$ is constrained to the following quasi-Newton condition:
\begin{align}\label{eq:quasi_newton_condition}
    H_{k+1} (\nabla f(x_{k+1}) - \nabla f(x_k)) = x_{k+1} - x_k,
\end{align}
which is inspired by Taylor series at point $x_{k+1}$
\begin{align*}
    \nabla f(x_k) - \nabla f(x_{k+1}) &= \nabla^2 f(x_{k+1}) (x_k - x_{k+1}) + o(\|x_k - x_{k+1}\|_2) \\
    x_{k+1} - x_k &\approx [\nabla^2 f(x_{k+1})]^{-1}(\nabla f(x_{k+1}) - \nabla f(x_k))
\end{align*}

There are several methods that use different rules to satisfy criteria \eqref{eq:quasi_newton_condition}~--- some of them (viz. DFP) are also presented in the method comparison table \ref{fig:fgm}. However, one of the most practically efficient quasi-Newton methods is BFGS (the results of which are also presented in the table \ref{fig:fgm}): 
\begin{gather*}
    x_{k+1} = x_k - h_k \cdot H_k \nabla f(x_k),\text{ where }h_k = \arg \min_{h > 0} f(x_k - h \cdot H_k \nabla f(x_k)),\\
    H_{k+1} = H_k + \frac{H_k \gamma_k \delta_k^\top + \delta_k \gamma_k^\top H_k}{\langle H_k \gamma_k, \gamma_k \rangle} - \beta_k \frac{H_k \gamma_k \gamma_k^\top H_k}{\langle H_k \gamma_k, \gamma_k \rangle},\\
    \text{where  }\beta_k = 1 + \frac{\langle \gamma_k, \delta_k \rangle}{\langle H_k \gamma_k, \gamma_k \rangle}, \gamma_k = \nabla f(x_{k+1}) - \nabla f(x_k), \delta_k = x_{k+1} - x_k, H_0 = I.
\end{gather*}
One of its practically valuable features is the stability to computational and line-search inaccuracies (there were also tested several methods of line-search, the most effective and economical turned out to be the method of quadratic interpolations, tuned for a certain fixed number of iterations.
However, this method is unsuitable for large-scale problems due to the large amount of memory required to store the matrix $H_k$. Therefore, in practice, the method of recalculating the $H_k$ matrix using only $r$ vectors $\gamma_k$ and $\delta_k$ from the last iterations is often used \cite{nocedal1980updating}, in this case $H_{k-r}$ is assumed to be equal to $I$. The described principle underlies the L-BFGS($r$) methods class with a $r$ memory depth.

Theoretically, it is known about quasi-Newton methods that the global rate of their convergence in the case of smooth convex problems does not exceed the estimates obtained for the classical gradient method, and in the vicinity of the minimum the rate of convergence becomes superlinear \cite{dennis1977quasi}. At the same time, in practical terms, the L-BFGS method (and its various modifications) is one of the most universal and effective methods of convex and even unimodal optimization \cite{liu1989limited,skajaa2010limited}, which allows us to assume the possibility of its effective application to the problem under consideration. Now, let us assess the practical efficiency of used version of the L-BFGS method: 
the results of work at various settings $r = 3, ..., 900$ are presented in table \ref{fig:fgm}. In addition to quasi-Newton methods, the table \ref{fig:fgm} presents the experimental results for various versions of Polyak method (Polyak) \cite{polyak1969}, Barzilai--Borwein method (BB) \cite{barzilai1988}, conjugate gradient method (CG) \cite{andrei40cg2008} and steepest descent method with zeroed small gradient components (Raider). 

As you can see from the presented data, the L-BFGS method actually demonstrates better performance compared to other methods. One of the unexpected results of this experiments is the special efficiency of the L-BFGS method in the case of using a large amount of information from past iterations. Classically, limited memory variants of the BFGS method have small optimal values of the history size, and are not so dependent on it, however, in this case, the best convergence rate of the L-BFGS method is achieved for value $r=900$, and with a further increase in this parameter, the result does not improve. This can be thought as one of the special and remarkable properties of the particular problem under consideration. 

\vspace{-0.7cm}
\hspace{-1.2cm}
\begin{minipage}[t]{\dimexpr.6\textwidth-.5\columnsep}
\raggedright

\begin{table}[H]
\centering
\caption{Full-gradient methods convergence, no time limit, residual model}
\label{fig:fgm}
\begin{tabular}{l||r|r|r|r}
\multirow{2}{*}{Method} & \multicolumn{4}{c}{Time to reach dB, sec.} \\ \cline{2-5}
& -30 dB & -35 dB & -37 dB & -39 dB \\ \hline \hline

{DFP(100)} & 11.02 & 50.22 & 83.46 & 2474.27 \\
{DFP(200)} & 10.27 & 41.48 & 79.75 & 1750.35 \\
{DFP(300)} & 10.24 & 41.32 & 69.49 & 2120.14 \\
{DFP(400)} & 10.29 & 42.07 & 70.27 & 1344.80 \\
{DFP(inf)} & 11.10 & 43.96 & 72.56 & 944.52 \\ \hline

{BFGS(100)} & 7.41 & 34.77 & \ok{30}{53.89} & \ok{10}{695.28} \\
{BFGS(200)} & 10.01 & 44.97 & 228.66 & 3747.02 \\
{BFGS(300)} & 10.70 & 48.00 & 202.40 & 3870.83 \\
{BFGS(400)} & 10.02 & 44.95 & 211.28 & 4104.48 \\
{BFGS(inf)} & 10.65 & 47.96 & 188.55 & \\ \hline

{LBFGS(3)} & \ok{30}{4.75} & \ok{30}{22.47} & 57.34 & \\
{LBFGS(10)} & \ok{60}{4.31} & \ok{60}{19.45} & \ok{30}{48.82} & 777.82 \\
{LBFGS(100)} & \ok{60}{4.04} & \ok{100}{16.81} & \ok{60}{36.09} & \ok{30}{513.34} \\
{LBFGS(300)} & \ok{60}{4.02} & \ok{100}{16.77} & \ok{100}{34.51} & \ok{60}{449.31} \\
{LBFGS(500)} & \ok{60}{4.01} & \ok{100}{16.80} & \ok{100}{34.63} & \ok{60}{424.39} \\
{LBFGS(700)} & \ok{60}{4.00} & \ok{100}{16.77} & \ok{100}{34.62} & \ok{100}{410.86} \\
{LBFGS(900)} & \ok{60}{4.01} & \ok{100}{16.75} & \ok{100}{34.49} & \ok{100}{399.52} \\

\end{tabular}
\end{table}

\end{minipage}
\begin{minipage}[t]{\dimexpr.6\textwidth-.5\columnsep}
\raggedleft
\vspace{1.0cm}
\begin{table}[H]
\centering
\begin{tabular}{l||r|r|r|r}
\multirow{2}{*}{Method} & \multicolumn{4}{c}{Time to reach dB, sec.} \\ \cline{2-5}
& -30 dB & -35 dB & -37 dB & -39 dB \\ \hline \hline

{SDM} & 25.83 & 925.92 & 5604.51 & \\ \hline

{Polyak(orig.)} & 7.39 & 123.93 & 345.10 & \\
{Polyak(v1)} & 24.41 & 1317.58 & & \\
{Polyak(v2)} & 19.09 & 782.29 & & \\ \hline

{BB(v1)} & 9.82 & 148.09 & 386.93 & \\
{BB(v2)} & 10.25 & 186.78 & & \\ \hline

{Raider(0.1)} & 40.15 & 2602.57 & & \\
{Raider(0.2)} & 482.80 & & & \\
{Raider(0.3)} & 266.83 & & & \\ \hline

{CG(HS)} & \ok{30}{4.50} & 27.31 & 61.94 & \\
{CG(FR)} & 15.31 & 55.44 & 109.25 & 4209.56 \\
{CG(PRP)} & \ok{100}{3.88} & \ok{30}{26.13} & 60.10 & \\
{CG(PRP+)} & \ok{60}{4.16} & \ok{30}{24.91} & 60.29 & \\
{CG(CD)} & \ok{60}{4.29} & \ok{30}{25.15} & 59.90 & \\
{CG(LS)} & \ok{60}{4.17} & \ok{30}{26.58} & 60.86 & \\
{CG(DY)} & 12.38 & 47.21 & 80.74 & 3721.43 \\
{CG(Nesterov)} & 28.86 & 114.86 & 277.88 & \\

\end{tabular}
\end{table}

\end{minipage}
\vspace{0.3cm}

\subsection{Flexible Gauss--Newton Method}

Let us now proceed another possible approach to solving the described problem, using the ideas underlying the Gauss--Newton method for solving the nonlinear least squares problem. The approach described in this section was proposed by Yu.E.~Nesterov in work\footnote{This paper is in print. The result of Nesterov's paper and our paper make up the core of the joint Huawei project.
The described below Method of Three Squares \cite{nesterov2021flexible} was developed as an attempt to beat L-BFGS (see fig.~\ref{fig:ssm_conv}). We repeat in this paper the main results of \cite{nesterov2021flexible} since they were developed for considered problem formulation and for the moment there is no possibility to read about these results somewhere else. Note, that recently some  results of the paper \cite{nesterov2021flexible} were generalized \cite{yudin2021flexible}. In particular, in \cite{yudin2021flexible} one can find more information about the Method of Three Squares.} \cite{nesterov2021flexible}. Let us reformulate the original problem \eqref{main_problem}. Consider a mapping $F: \mathbb{R}^n \rightarrow \mathbb{R}^m $ of the following form: 
\[
 F(x) := (F_1(x), \dots, F_m(x)),
\]
where each component represents the discrepancy between the approximation obtained by the model and the exact solution for each of the objects of the training set: $F_i(x) := [\mathcal{M}_\theta (x)]_i - \overline{y}_i$. Then the original problem can be reduced to solving the following least squares problem: 
\begin{equation} \label{prob-Main}
    \min_{x \in \R^n} \{f_1(x) := \| F(x) \|_2 \}.
\end{equation}
We additionally require only the Lipschitz smoothness of the functional $F$ (note that throughout the analysis of the method, the requirement of convexity will not be imposed on the functional, that is, the presented convergence estimates are valid in non-convex generality):
\[
\| F'(x) - F'(y) \|_2 \leq L_F \| x - y \|_2, \quad x, y \in \R^n,
\]
where $F'(x) = \left(\frac{\partial F_i(x)}{\partial x_j}\right)_{i, j}$ is a Jacobian. Under these assumptions, one can prove the following lemma on the majorant for the initial function $f_1$: 
\begin{lemma}{\cite{nesterov2021flexible}}
Let $x$ and $y$ be some points from $\R^n$, $L \geq L_F$,and $f_1(x) > 0$. Then
\begin{equation} \label{major1}
 f_1(y) \leq \hat \psi_{x,L}(y) := \frac{1}{2 f_1(x)} \Big[ f_1^2(x) + \| F(x) + F'(x)(y-x) \|_2^2 \Big] + \frac{L}{2} \| y - x \|_2^2.
\end{equation}
\end{lemma}

Let us assume for a moment that we know an upper bound $L$ for the Lipschitz constant~$L_F$. Then the last inequality in \eqref{major1} leads to the following method:

\begin{align}\label{met-3S}
&\textbf{\hspace{25ex} Method of Three Squares (\cite{nesterov2021flexible})}\\
&\nonumber x_{k+1} = \arg\min\limits_{y \in \R^n} \Big\{ {\frac{L}{2}}
\| y - x_k \|_2^2 + {\frac{1}{2} f_1(x_k)} \Big[ f_1^2(x_k) +
\| F(x_k) + F'(x_k)(y-x_k) \|_2^2 \Big] \Big\}
\end{align}

\bigskip
\noindent
The global convergence of this method is characterized by the following theorem. 

\begin{theorem}{\cite{nesterov2021flexible}} \label{th-Rate3S}
Let us assume that the function $F(\cdot)$ is uniformly
non-degenerate: $F'(x) F'(x)^\top \succeq \mu I_m$ for all $x \in {\cal F}_0 := \{ x \in \R^n: \; f_1(x) \leq f_1(x_0) \}$. If in the method (\ref{met-3S}), we choose $L \geq L_F$, then it converges linearly to the solution of equation $F(x)=0$:
\[
 f_1(x_k) \leq f_1(x_0) \cdot \exp\left\{ - \frac{\mu k}{2(L f_1(x_0) + \mu)} \right\}, \quad k \geq 0.
\]
At the same time, for any $k \geq 0$ we have
\[
 f_1(x_{k+1}) \leq \frac{1}{2} f_1(x_k) + \frac{L}{2 \mu} f^2(x_k).
\]
Thus, the coefficient for asymptotic local linear rate of
convergence for this method is $\frac{1}{2}$.
\end{theorem}

\bigskip
\noindent
If we relax the assumptions of theorem \ref{th-Rate3S},
then we can estimate the rate of convergence of this
method to a stationary point of problem (\ref{prob-Main}).
Denote
$$
f_2(x) := f_1^2(x) = \| F(x) \|_2^2.
$$
\begin{theorem}{\cite{nesterov2021flexible}} \label{th-Stat}
Suppose that the function $F$ has uniformly bounded derivative: $\|F'(x)\|_2 \leq M_F$ for all $x \in {\cal F}_0$. If in the method (\ref{met-3S}) $L \geq L_f$, then for any $k \geq 0$ we have
\[
 f_2(x_k) - f_2(x_{k+1}) \geq \frac{1}{8( L f_1(x_0) + M_F^2)} \| \nabla f_2(x_k) \|_2^2.
\]
\end{theorem}

\noindent
Thus, under very mild assumption (bounded derivative), we can prove that the measure of non-stationarity $\| \nabla f_2(\cdot) \|_2^2$ is decreasing as follows \cite{nesterov2021flexible}:
\[
 \min_{0 \leq i \leq k} \|\nabla f_2(x_k) \|_2^2 \leq \frac{8 f_2(x_0) (L f_1(x_0) + M_F^2)}{k+1}, \quad k \geq 0.
\]

We will also consider the following enhanced version of
method (\ref{met-3S}).

\begin{align}\label{met-NGN}
&\textbf{\hspace{13ex} Non-Smooth Gauss-Newton Method (\cite{nesterov2021flexible})}\\
&\nonumber x_{k+1} = \arg\min\limits_{y \in \R^n} \Big\{
\psi_{x_k,L}(y) \; := \; {\frac{L}{2}} \| y - x_k \|_2^2 + \|
F(x_k) + F'(x_k)(y-x_k) \|_2 \Big\}
\end{align}

\bigskip
\noindent
Let us describe its convergence properties.
\begin{theorem}{\cite{nesterov2021flexible}} \label{th-NGN}
Let us choose in the method (\ref{met-NGN}) $L \geq L_f$.
\begin{enumerate}
\item If function $F$ is uniformly non-degenerate: $F'(x) F'(x)^\top \succeq \mu I_m, \; x \in {\cal F}_0$, then method (\ref{met-NGN}) converges linearly to the
solution of equation $F(x)=0$:
\[
 f_1(x_k) \leq f_1(x_0) \cdot \exp\left\{ - \frac{ \mu k}{2(L f_1(x_0) + \mu)} \right\}, \quad k \geq 0.
\]

At the same time, it has local quadratic convergence:
\[
 f_1(x_{k+1}) \leq \frac{L}{2 \mu} f_1^2(x_k), \quad k \geq 0.
\]

\item Suppose that function $F$ has uniformly bounded
derivative: $\|F'(x)\|_2 \leq M_F, \; x \in {\cal F}_0$, then for any $k \geq 0$ we have
\[
 f_2(x_k) - f_2(x_{k+1}) \geq \frac{1}{8( L f_1(x_0) + M_F^2)} \|\nabla f_2(x_k)\|_2^2.
\]
\end{enumerate} 
\end{theorem}

Now, let us consider the normalized versions of introduced objective functions:
$
\hat f_1(x) := \frac{1}{\sqrt{m}} f_1(x) = \left[\frac{1}{m} f_2(x) \right]^{1/2},\; \hat f_2(x) := \frac{1}{m} f_2(x).
$
This normalization allows us to consider $m \to \infty$. 
Moreover, the objective function in this form admits stochastic approximation. Therefore, let us describe a stochastic variant of method (\ref{met-3S}).

\begin{align}\label{met-3SS}
&\textbf{\hspace{16ex} Method of Stochastic Squares (\cite{nesterov2021flexible})}\\
&\nonumber \mbox{a) Choose $L_0 > 0$ and fix the batch size $p \in \{0,\dots,m\}$.}\\
&\nonumber \mbox{b) Form $I_k \subseteq \{1, \dots, m\}$ with $|I_k| = p$ and define $G_k := \{F'_i(x_k), i \in I_k \}$.}\\
&\nonumber \mbox{c) Define $\vf_k(y) := \hat f_1(x_k) + \la \hat f'_1(x_k), y - x_k \ra + {\frac{1}{2 \hat f_1(x_k)}} \left({\frac{1}{p}} \| G^T_k (y - x_k) \|_2^2 \right)$.}\\
&\nonumber \mbox{d) Find the smallest $i_k \geq 0$ such that for the point}\\
&\nonumber \hspace{12ex} T_{i_k} = \arg\min\limits_y \left\{ \psi_{i_k}(y) := \vf_k(y) + {\frac{2^{i_k} L_k}{2}} \| y - x_k\|_2^2 \right\}\\
&\nonumber \text{ we have }\hat f_1(T_{i_k}) \leq \psi_{i_k}(T_{i_k})\\
&\nonumber \mbox{e) Set $x_{k+1} = T_{i_k}$ and $L_{k+1} = 2^{i_k-1} L_k$.}
\end{align}

\subsubsection{Numerical experiments}

All variants of the Gauss--Newton method described above were implemented taking into account the above remarks about the possibilities of their effective implementation. In a series of numerical experiments, it was tested the practical efficiency of two described full-gradient methods \eqref{met-3S} (3SM) and \eqref{met-NGN} (NsGNM), and the stochastic method \eqref{met-3SS} (SSM) for various batch sizes $p$. The results are presented in the table \ref{3sm_compare}, along with the results of the Gauss--Newton method in the Levenberg--Marquardt version (LM) \cite{marquardt1963algorithm}. It can be seen from the presented results that the Three Squares method demonstrates better performance than the Levenberg--Marquardt method, and, moreover, batching technique significantly accelerates the convergence of the proposed scheme. The best setting of the method with\footnote{Note, that $p\sim n$ can be easily explained by the following observation. In this regime Jacobian calculation $\sim pn^2$ has the same complexity as Jacobian inversion $\sim n^3$. It means that there is no reason to choose $p$ large, but $p \ll n$. If $p$ is large we can consider $p$ to be greater than $n$ since the complexity of each iteration include $n^3$ term anyway.} $p = 6n$ demonstrates a result that exceeds the performance of the L-BFGS method, starting from the mark of $-39$ dB of the quality metric (see fig. \ref{fig:ssm_conv}).

\vspace{-0.5cm}
\hspace{-1.0cm}
\begin{minipage}[t]{\dimexpr.6\textwidth-.5\columnsep}
\raggedright

\begin{table}[H]
\centering
\caption{Gauss-Newton methods convergence, residual model}
\label{3sm_compare}
\begin{tabular}{l||r|r|r|r}
\multirow{2}{*}{Method} & \multicolumn{4}{c}{Time to reach dB level, sec.} \\ \cline{2-5}
& -30 dB & -35 dB & -37 dB & -39 dB \\ \hline \hline

{LM(1)} & 919.28 & 1665.08 & 3165.61 & 14270.56 \\
{LM(2)} & 760.54 & 2277.11 & 5534.27 & \\
{LM(3)} & 550.06 & 1494.13 & 3584.74 & \\ \hline

{3SM} & 633.88 & 1586.69 & 1747.28 & 4616.80 \\
{NsGNM} & 762.99 & 1626.18 & 2265.74 & 9172.11 \\ \hline

{SSM($0.3 n$)} & 93.57 & 271.63 & 864.04 & \\
{SSM($0.6 n$)} & 50.43 & 145.60 & 437.22 & \\
{SSM($1 n$)} & 48.79 & 104.31 & 186.37 & 12469.39 \\
{SSM($2 n$)} & 60.43 & 91.62 & 141.99 & 4407.34 \\
{SSM($3 n$)} & 37.13 & \ok{60}{77.75} & \ok{30}{116.78} & 1118.43 \\
{SSM($4 n$)} & 53.71 & 96.04 & \ok{60}{110.08} & 570.92 \\
{SSM($5 n$)} & 42.65 & \ok{30}{88.14} & \ok{60}{110.00} & \ok{100}{315.28} \\
{SSM($6 n$)} & 47.32 & \ok{100}{72.30} & \ok{100}{89.08} & \ok{100}{314.95} \\
{SSM($7 n$)} & 44.16 & 109.33 & 137.58 & 382.58 \\
{SSM($8 n$)} & 39.71 & 92.33 & \ok{30}{114.27} & 400.38 \\
{SSM($10 n$)} & 46.65 & 120.24 & 145.04 & \ok{60}{331.06} \\
{SSM($11 n$)} & 48.43 & 111.66 & 150.10 & 343.44 \\
\end{tabular}
\end{table}

\end{minipage}
\begin{minipage}[t]{\dimexpr.6\textwidth-.5\columnsep}
\raggedleft

\begin{figure}[H]
\centering
\includegraphics[width=\linewidth]{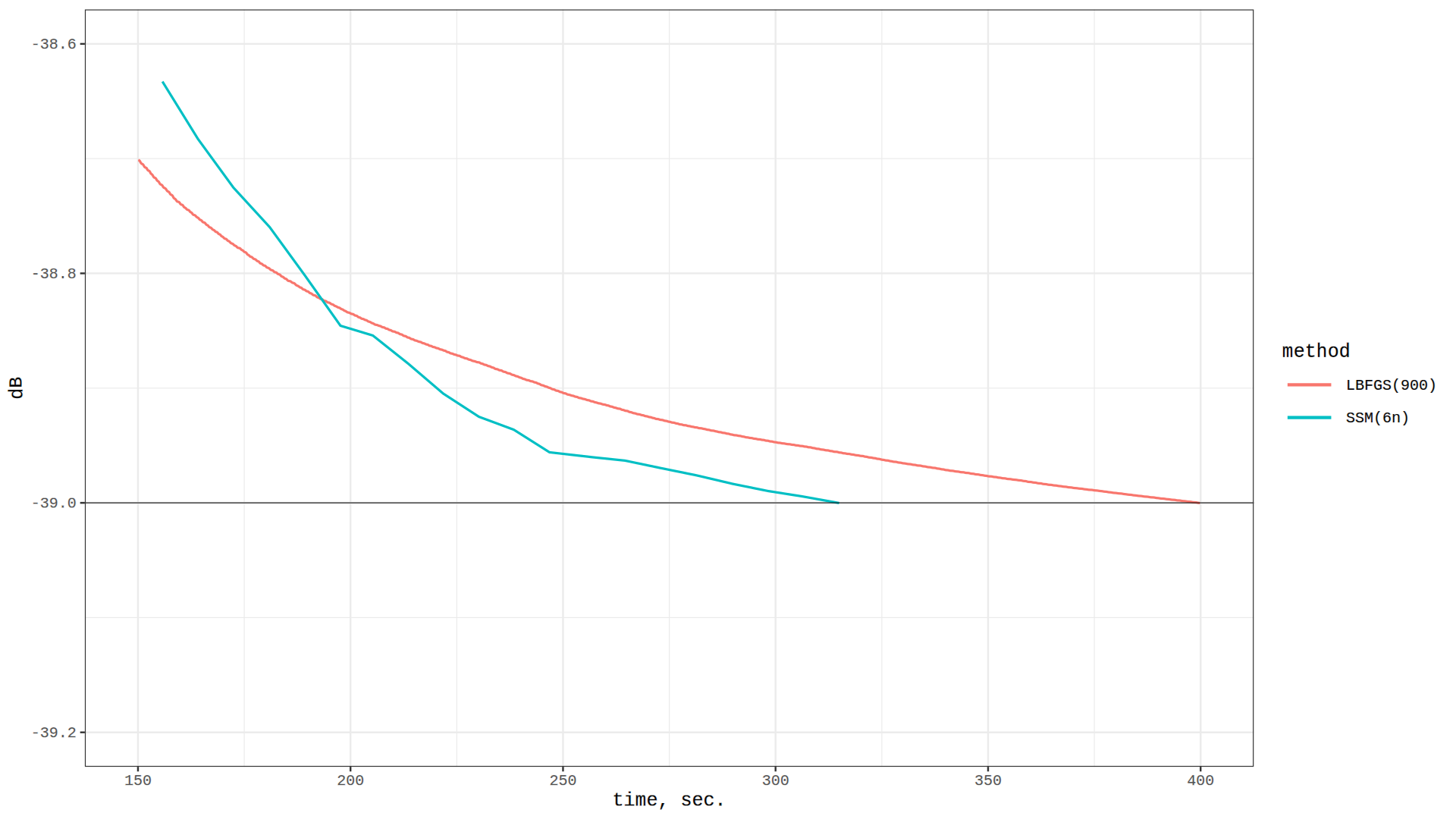}
\caption{Convergence of Stochastic Squares Method (SSM) and L-BFGS method, residual model}
\label{fig:ssm_conv}
\end{figure}

\end{minipage}

\subsection{Stochastic methods}

Let us now turn to a class of methods that are especially widely used for problems related to training models represented in the form of large computational graphs (in particular, neural networks)~--- stochastic gradient methods. In addition to the repeatedly confirmed practical efficiency, the motivation for applying stochastic methods to the existing problem is a significant saving of time when evaluating the function for only one of the terms of the sum-type functional at our disposal. Indeed, consider the calculation complexity for the one term~--- it requires not more than $4 \cdot (R_y+1)\cdot(R_y N_z + 2 M_z P_z)\simeq 4\cdot 10^3$ arithmetical operations, whereas the complexity for the full sum is $4 \cdot 2m \cdot(R_y N_z + 2 M_z P_z) \simeq 3 \cdot 10^8$ a.o., that almost in $\sim m$ times more. Further, according to theory of automatic differentiation, for the particular computational graph the calculation of gradient is not more than 4 times more expensive than calculation of the function value \cite{griewank1989automatic,nocedal2006numerical}, although it is necessary to store the entire computational graph in RAM. 


In our experiments, we apply stochastic methods to the problem under consideration. Along with classical stochastic gradient method (SGD), there were tested various modifications of adaptive and adaptive momentum stochastic methods (Adam, Adagrad, Adadelta, Adamax). Adaptability of these methods lies in the absence of the need to know the smoothness constants of the objective function, which is especially effective in deep learning problems \cite{zhang2019adam}, although theoretically has no advantages in terms of convergence rate. Moreover, a number of variance reduction methods were applied (SVRG, SpiderBoost). SpiderBoost matches theoretical lower bound for the number of iterations sufficient to achieve a given accuracy in the class of stochastic methods under the assumption that all terms are smooth. Note, that the code of most of these methods is free available at GitHub: \url{https://github.com/jettify/pytorch-optimizer}.

As one can see from the figures below, adaptive stochastic methods, in particular Adam, show the most effective (among stochastic methods) convergence for the considered model (see fig.~\ref{fig:adaptive}). At the same time, the use of variance reduction methods does not allow achieving any acceleration of the convergence (see fig.~\ref{fig:svrg}). Note also that variance reduction methods are inferior in convergence rate to the standard SGD method also in terms of the number of passes through the dataset.

Experiments also show that the efficiency of stochastic algorithms (in comparison with full-gradient methods) significantly depends on the dimension of the parameter space of the model used. Moreover, for models with a large number of blocks, the rate of convergence of Adam-type algorithms is slower than for methods of the L-BFGS type, due to the slowdown in convergence with an increase in the number of iterations of the method. It is important to note, however, that losing in the considered setting in terms of the depth and rate characteristics of convergence, stochastic methods show the advantage of being more resistant to overfitting \cite{amir2021sgd} (due to their randomized nature). At the same time, the stochastic methods in the current version are especially valuable for the possibility of using them for online training of the model. Indeed, the main application of the solution to the problem posed at the beginning of the article is to optimize the DPD function, however, it is quite natural that with a change in the characteristics of the input signal over time, the optimal parametrization of the model can smoothly change, so it is necessary to adjust the model to the new data. Modification of full-gradient methods for their efficient operation on mini-batches is a separate non-trivial problem, while stochastic methods provide such opportunities out of the box. Also, stochastic methods are more convenient for their hardware implementation, since they do not require storing long signal segments for training.

\vspace{-0.6cm}
\hspace{-1.0cm}
\begin{minipage}[t]{\dimexpr.58\textwidth-.5\columnsep}
\raggedright

\begin{table}[H]
\centering
\caption{Stochastic methods convergence, residual model }
\begin{tabular}{l|rr||r|r}
\multirow{2}{*}{Method} & \multicolumn{2}{c||}{\multirow{2}{*}{Setup}}
& \multicolumn{2}{|c}{dB} \\
\cline{4-5}
& & & $t=0$ sec. & $t=300$ sec. \\
\hline
ASGD & (128,& 10.0) & \multirow{8}{*}{-15.616} & -28.976 \\
Adadelta & (2048,& 10.0) & & -32.697 \\
Adagrad & (2048,& 0.01) & & \ok{30}{-36.608} \\
Adam &(2048,& 0.001) & & \ok{100}{-38.129} \\
Adamax& (2048,& 0.01) & & \ok{60}{-37.934} \\
RMSprop& (2048, &0.001) & & -36.499 \\
SGD & (128,& 10.0) & & -34.061 \\
FastAdaptive & & & & -36.273 \\
\end{tabular}
\end{table}

\end{minipage}
\hspace{0.6cm}
\begin{minipage}[t]{\dimexpr.45\textwidth-.5\columnsep}
\raggedleft

\begin{figure}[H]
\centering
\includegraphics[width=\linewidth]{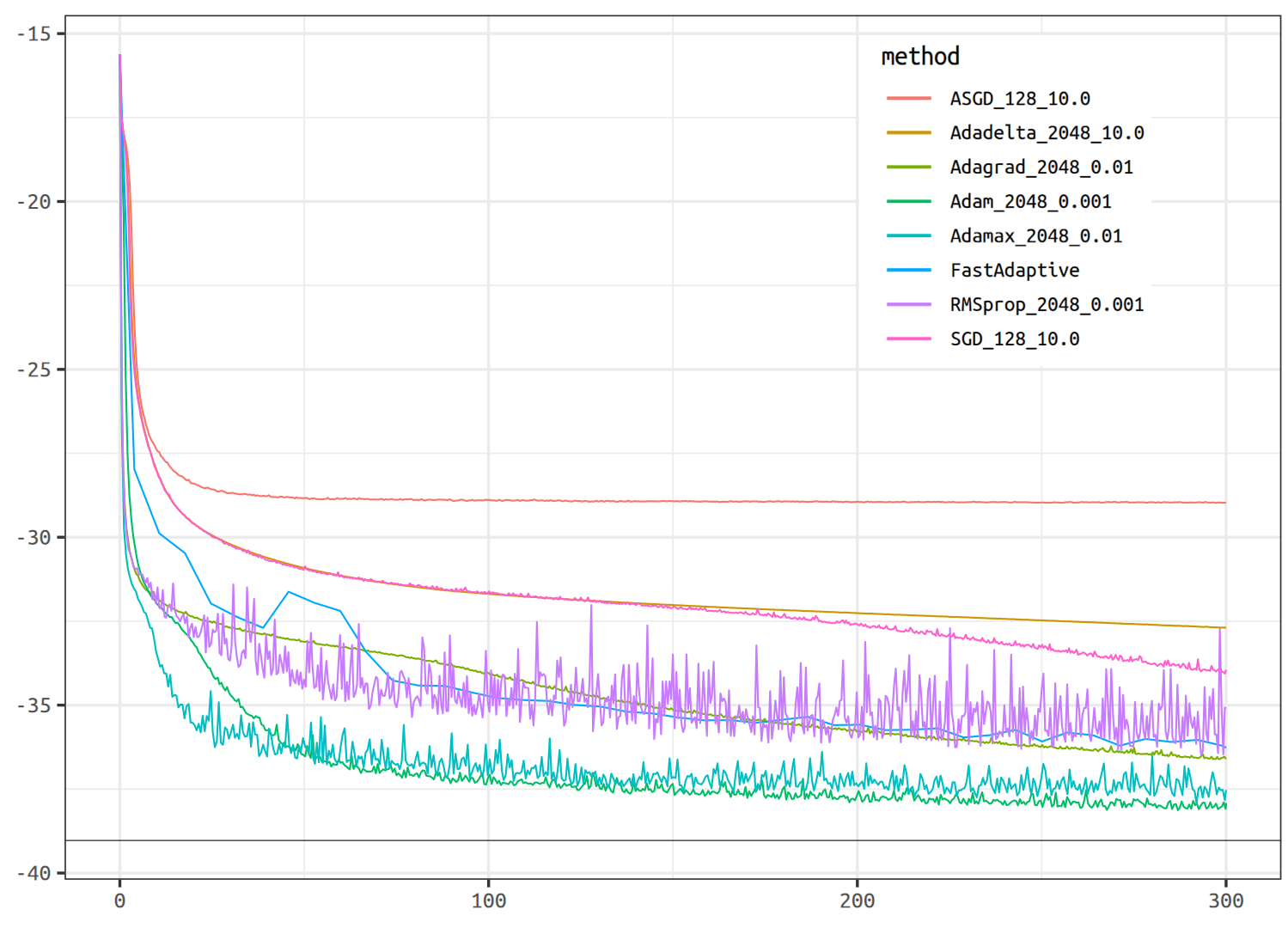}
\caption{Stochastic methods convergence, residual model}
\label{fig:adaptive}
\end{figure}

\end{minipage}
\vspace{-0.6cm}

\begin{minipage}[t]{\dimexpr.47\textwidth-.5\columnsep}
\raggedright
\begin{figure}[H]
    \centering
    \vspace{-0.6cm}
    \includegraphics[width=\linewidth]{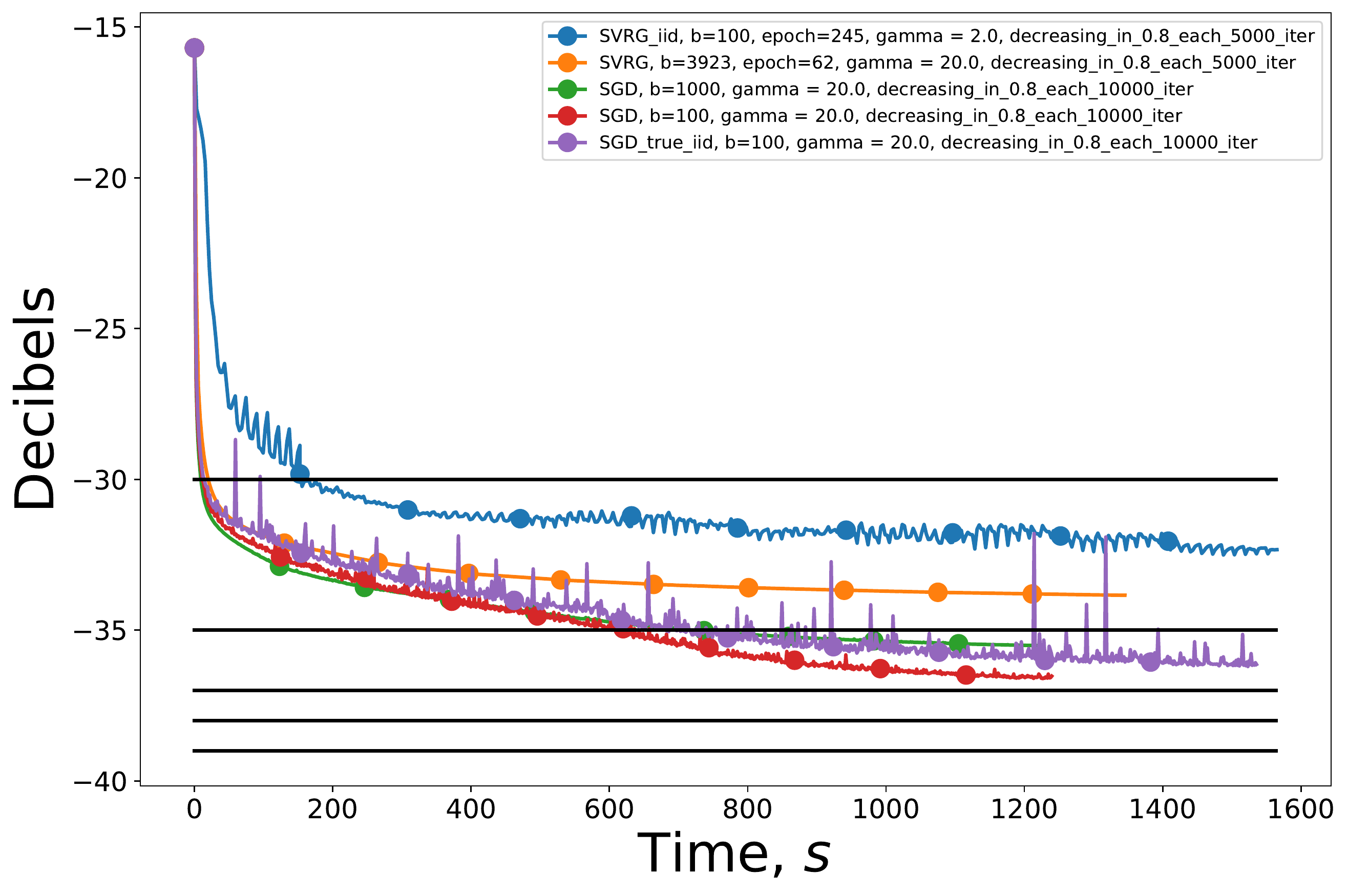}
    \caption{Comparison of various sampling strategies for SVRG and SGD in terms of running time (in seconds) to reach the predefined threshold ($-30,-35,-37,-38$ dB).}
    \label{fig:svrg}
\end{figure}

\end{minipage}
\begin{minipage}[t]{\dimexpr.2\textwidth}

\end{minipage}
\begin{minipage}[t]{\dimexpr.47\textwidth-.5\columnsep}
\raggedleft

\begin{figure}[H]
    \centering
    \vspace{-0.6cm}
    \includegraphics[width=\linewidth]{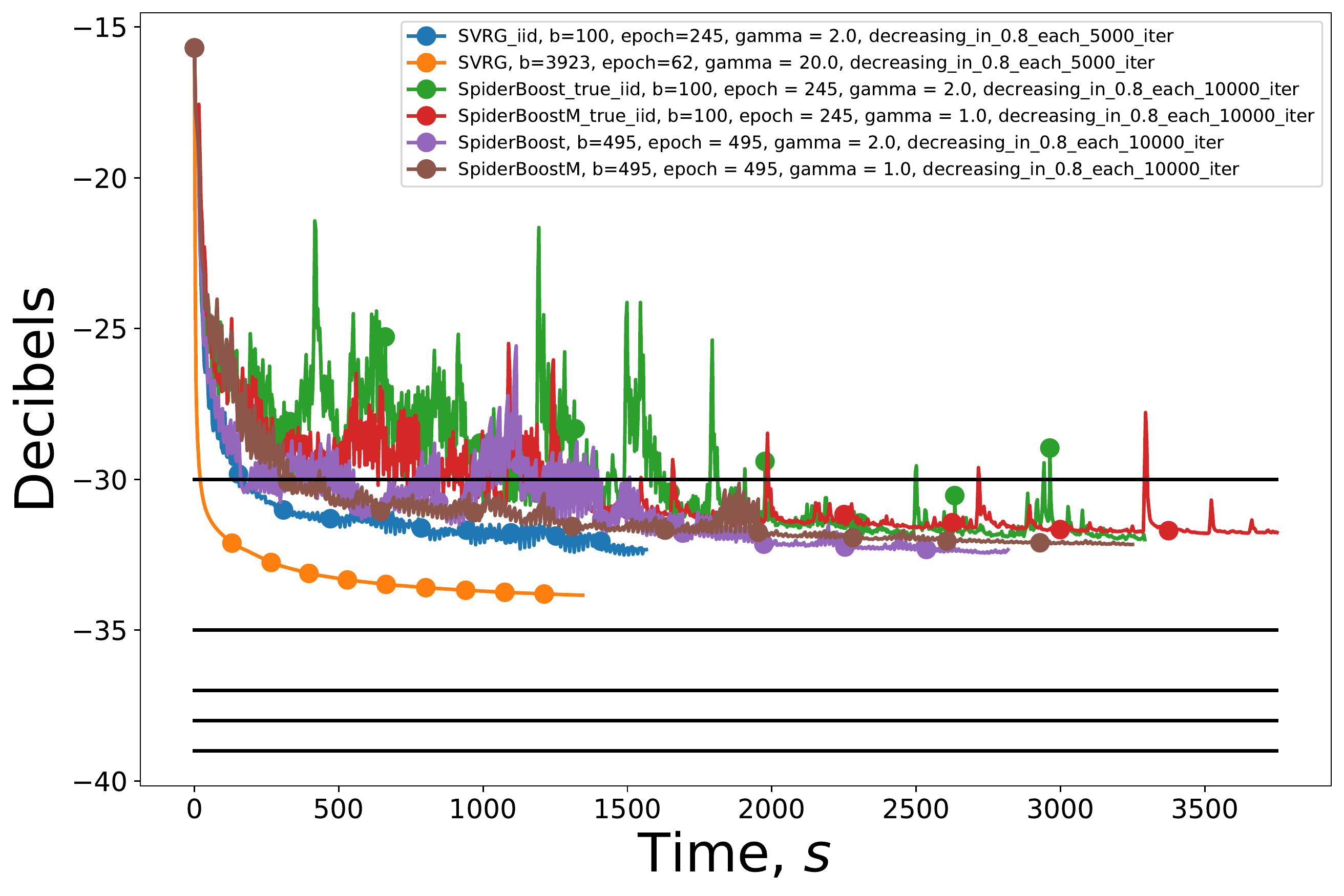}
    \caption{Comparison of various sampling strategies for SVRG and SpiderBoost in terms of running time (in seconds) to reach the predefined threshold ($-30,-35,-37,-38$ dB).}
\end{figure}

\end{minipage}
\vspace{-0.4cm}
\subsection{Global algorithms}

\subsubsection{Multi-start tests} \label{multistart}
\hfill

From the point of view of characterizing the model under consideration, it is useful to investigate its global characteristics. In this section, we consider the issue of optimization methods operation stability with a random choice of the starting point for them, for two variants of the model: a direct model and a model utilizing skip connections. We perform numerous multi-start tests for both original and residual model according to the following principle: firstly, a set of start points are uniformly drawn from n-orthotope $[-0.1, 0.1]^n$ and then L-BFGS(900) method works starting from every point until it fails to find relaxation for the next step (neither the maximum number of iterations nor the maximum working time nor any other stop criteria is specified). 

The performed numerical experiments showed that the original model has a significant instability of the results depending on the starting point. The difference between the worst and best obtained solutions exceeds $5$ dB and the effect of the multi-start scheme is quite noticeable (see fig.~\ref{fig:multi-original}). Residual model on the contrary showed high enough stability~--- for the majority of starting points the result of optimization gives values below $-39$ dB. The difference between the worst and best solutions is about $2-2.5$ dB and the effect of the multi-start scheme is not too significant (see fig.~\ref{fig:multi-residual}).

\vspace{-0.8cm}
\hspace{-1.0cm}
\begin{minipage}[t]{\dimexpr.5\textwidth-.5\columnsep}
\raggedright

\begin{figure}[H]
    \centering
    \includegraphics[width=0.95\linewidth]{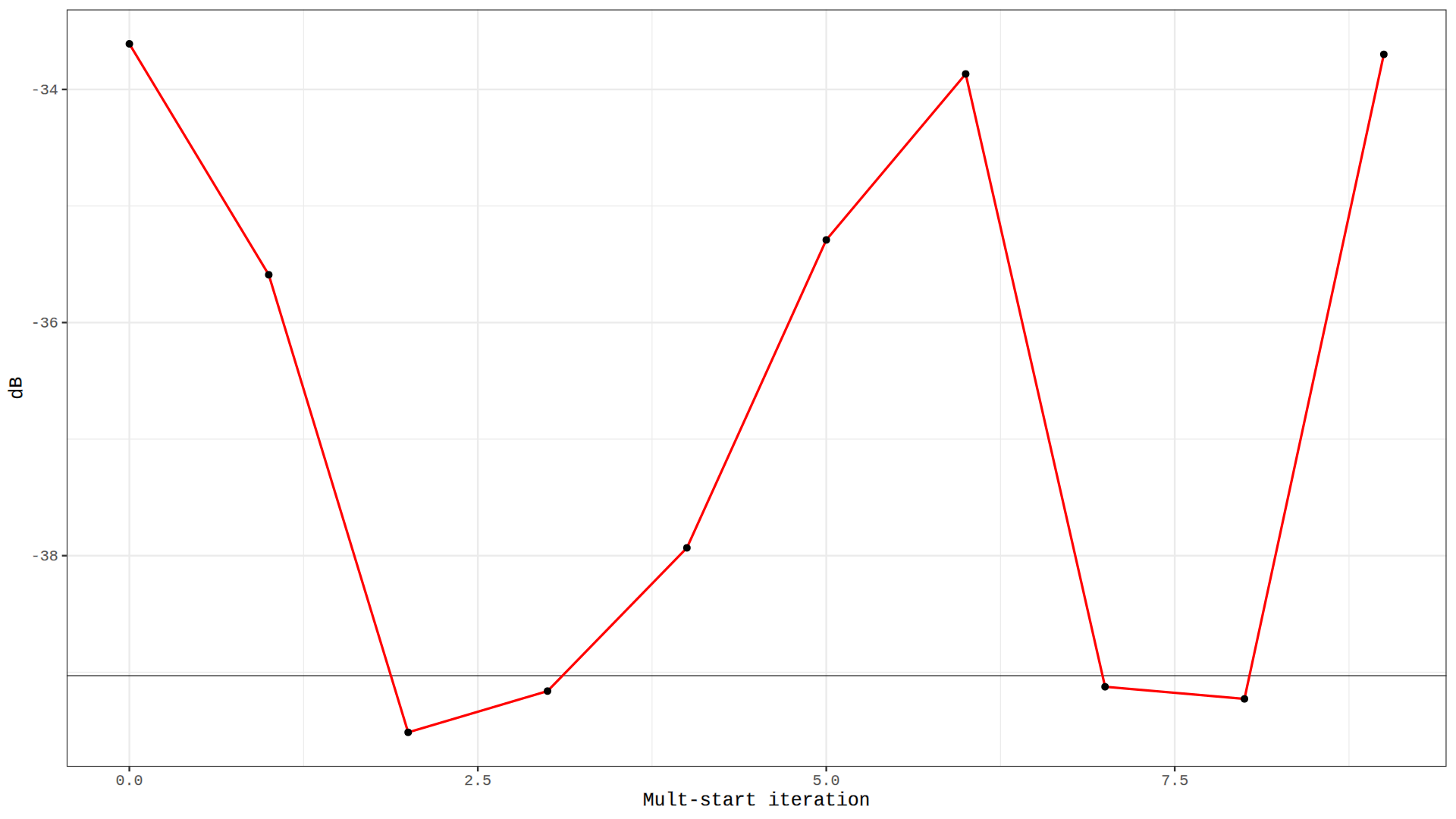}
    \caption{Multi-start results, original model}
    \label{fig:multi-original}
\end{figure}

\end{minipage}
\hspace{0.4cm}
\begin{minipage}[t]{\dimexpr.5\textwidth-.5\columnsep}
\raggedleft

\begin{figure}[H]
    \centering
    \includegraphics[width=0.95\linewidth]{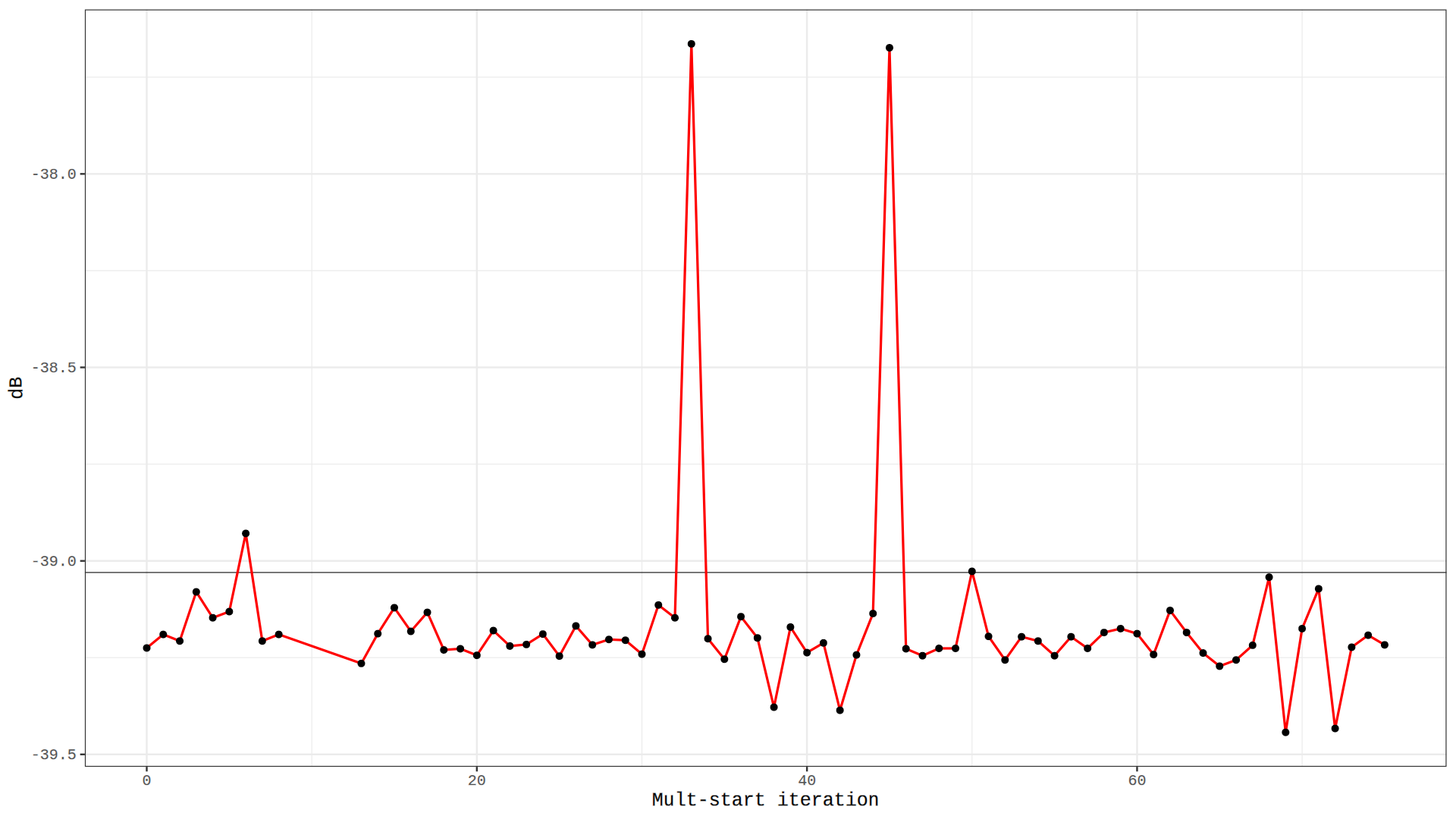}
    \caption{Multi-start results, residual model}
    \label{fig:multi-residual}
\end{figure}

\end{minipage}
\vspace{-0.4cm}

\subsubsection{Simulated annealing}
\hfill\\
Optimization method considered in this division is an adaptation of the Metropolis--Hastings algorithm to the analogy of thermodynamic system states evolution. The simulated annealing method is a non-local optimization method that allows one to construct, according to the same principle, many different modifications of the method for various particular classes of extremal problems, both continuous and discrete. The variant of the method used in this article is based on a particular Metropolis--Kirkpatrick procedure \cite{kirkpatrick1983optimization}. To date, there are quite a few modifications of algorithms such as simulated annealing, for which can be formulated some strict theoretical statements on their convergence properties \cite{hajek1988cooling,nolte2000note}. This method is used for an extremely large variety of practical problems, including some NP-hard problems and significantly global optimization problems, where it turns out to be quite effective. 

There were implemented two variants of the considered scheme~--- the classic one, using only a value of the function (v1), and the modified one, in which the optimization trajectory is built based on the conjugate gradient method (v2). The decision to conduct a jump in the trajectory, which allows tunneling from one local solution to another, is made either by the built-in cooling strategy or forcibly through a priori specified number of iterations. 
Numerical experiments have shown that the value of the jump lengths is, in a sense, a parameter that regulates the globality of the method: at small values of the parameter, it is possible to achieve a faster local convergence of the method, while at large values, the chance of jumping into the region of attraction of a neighboring, possibly better, local minimum increases. It turned out that in the v2 version of the method, the number of local descent iterations must be related to the length of the jumps: the longer jump length is set, the greater the number of local descent iterations must be performed in order to preserve the relaxation of iterations (\ref{alg:sa}). Note also that the quality of the method depends significantly on the choice of the starting point.

\subsection{Differential evolution}
\vspace{-0.5cm}
\begin{table}[H]
\centering
\caption{Non-convex (global) methods performance, residual model}
\label{sa_de}
\begin{tabular}{l||r|r|r|r}
\multirow{2}{*}{Method} & \multicolumn{4}{c}{Time to reach dB level, sec.} \\ \cline{2-5}
& -30 dB & -35 dB & -37 dB & -39 dB \\ \hline \hline

{SA(0.1, 1500)} & 80.53 & 80.53 & 80.53 & \ok{65}{161.04} \\ \hline
{SA(0.1, 2500)} & 136.47 & 136.47 & 136.47 & \ok{80}{136.47} \\ \hline
{DE(5, 1000)} & 658.71 & 658.71 & 658.71 & 2377.19 \\
{DE(50, 100)} & 655.39 & 1879.50 & 2418.66 & \\
{DE(10, 2000)} & 1107.28 & 1107.28 & 1107.28 & 1107.28 \\
{DE(100, 200)} & 996.40 & 996.40 & 2846.19 & 9192.72 \\
\end{tabular}
\end{table}
\vspace{-0.5cm}

Differential evolution is another non-local, bio-inspired, population type method that based on the addition of information recorded in four randomly selected individuals from the population (three ``mothers'' and one ``father''). This approach was proposed in \cite{storn1995differrential}, many details and variations of the algorithm can be found in the book \cite{price2006differential}.

From the results of the experiments, also presented in the summary table~\ref{sa_de}, it can be seen that the population size directly affects the convergence rate of the algorithm: large values increase the degree of biodiversity (in the original sense, the degree of globalization of the method), but increases the complexity of each iteration. This obstacle can be technically circumvented by using multiple GPU's, since the processing for each element of the population does not depend on the other members. The number of iterations of the local algorithm also affects the quality of the algorithm: too small values degrade the relaxation properties of the iteration of the method.(\ref{alg:diff_evo}) The starting point affects the quality of work rather poorly, since the method is initially designed for the presence of many different instances in the population.

}
\section{Model Tuning}
{

\subsection{Model Structure}

A naturally arising interesting question is the optimal configuration of the used Wiener--Hammerstein model, that is, the potential approximating capabilities of the model itself~--- this is useful both from a practical point of view (if the optimal scheme turns out to be significantly more efficient than an arbitrarily chosen one), and with an abstract (comparison may reveal new dependencies of the quality of the approximation on the structure of the model, characterizing the tendency to overfitting or the specifics of the data, which could be used in further studies).

Several experiments have been carried out using optimization algorithm L-BFGS, the results of which are presented in the table~\ref{tab:configuration}. They can be summarized by the following points (the notation of section \ref{model_subs} is used): 
\begin{enumerate}
  \item Layers number $L$ is much more important hyperparameter than $R$, so we should prefer to increase $L$ instead of $R$.
  
  \item If some model has $R > L$, then by simple swapping this values we can generally get better results with the same number of parameters.

  \item The original model settings ($L = 2$, $R_{z/y} = 11$) give the worst results. The model with $L = 8$, $R_{z/y/..} = 2$ gives better results with fewer parameters count ($1072$ vs $1474$).
\end{enumerate}
\vspace{-0.8cm}
\hspace{0.6cm}
\begin{minipage}[t]{\dimexpr.65\textwidth-.5\columnsep}
\raggedright

\begin{table}[H]
\centering
\caption{Optimization results for different network setups, residual model}
\label{tab:configuration}
\begin{tabular}{l|r|r||r|r|r}
\multirow{2}{*}{Layers} & 
\multirow{2}{*}{$R$} & Parameters, &
\multicolumn{3}{c}{dB} \\
& & (complex) & $t = 0$ sec. & $t = 300$ sec. & $t = t_{end}$  \\
\hline

3 & 8 & \multirow{4}{*}{1608} & -15.666 & -40.552 &
-41.087\\

4 & 6 & & -15.577 & -41.048 &
-42.029  \\

6 & 4 & & -15.742 & \ok{30}{-41.727} &
\ok{30}{-43.375}  \\

8 & 3 & & -15.781 & -41.209 &
\ok{100}{-43.790} \\
\hline

2 &
11 &
\multirow{2}{*}{1474} &
-15.477 &
\bad{-38.920} &
\bad{-39.141}\\

11 & 2 & & -15.693 & -40.835 &
\ok{60}{-43.707} \\
\hline

3 & 7 & \multirow{2}{*}{1407} & -15.549 & -40.429 & -40.949 \\

7 & 3 & & -15.865 & -40.986 &
-42.690\\
\end{tabular}
\end{table}

\end{minipage}
\hspace{0.4cm}
\begin{minipage}[t]{\dimexpr.1\textwidth-.5\columnsep}
\raggedleft

\vspace{1.15cm}
\begin{table}[H]
\centering
\begin{tabular}{|r}
\multirow{2}{*}{$t_{end}$, sec.} \\
\\
\hline

1258.257 \\

1195.800 \\

2029.769 \\

2511.681 \\
\hline

1405.771 \\

2783.392 \\
\hline

1191.327 \\

946.900 \\
\end{tabular}
\end{table}

\end{minipage}
\vspace{0.3cm}
\subsection{Initialization}

Computational graphs parameters optimization problems have a number of specific features that complicate the operation of the numerical methods applied to them. In the case of classical neural networks, corresponding training problems are characterized by significant multimodality and the presence of complex ravines \cite{auer1996exponentially,choromanska2015loss,garipov2018loss}. If we consider in more detail the used Wiener--Hammerstein model, we can assume that the corresponding functional has extensive plateaus (regions of functional's regularity)~--- this follows from the model's output structure of the ``sum of products'' form. This raises the question of the correct initialization of the model parameters, at which the initial point of the methods trajectories could lie near the attraction region of the potential global minimum, and not in an arbitrary region of the parametric space with possibly poor local problem's properties. It is worth noting, however, that many random initialization strategies would be inappropriate in this case. The fact is that deep computational graphs during training are also subject to the phenomena of gradient vanishing and gradient exploding, which significantly impairs the convergence of methods and the trainability of the model in general. In view of this, we would also like to have a method for initializing weights, which allowing to maintain the variance of the values at the output of each layer approximately equal to the variance of the values at its input. 

In this research we considered the initialization techniques which is used in classical neural networks: Xavier and He initialization (the corresponding formulas presented in a table below), which, however, are applicable as effective heuristics in the case of the model under consideration.

\begin{table}[H]
\centering
\vspace{-0.5cm}
\begin{tabular}{c c}
\textbf{Xavier initialization}\vspace{0.15cm} & \quad\textbf{Initialization He}\vspace{0.15cm} \\ 

$\displaystyle w_{i} \sim U\left[-\frac{\sqrt{6}}{\sqrt{n_{\mathrm{in}}+n_{\mathrm{out}}}}, \frac{\sqrt{6}}{\sqrt{n_{\mathrm{in}}+n_{\mathrm{out}}}}\right]$ & $\quad\displaystyle w_{i} \sim \mathcal{N}\left(0, \sqrt{\frac{2}{n_{\mathrm{in}}}}\right)$

\vspace{-0.5cm}
\end{tabular}
\end{table}

These techniques rely on some theoretical analysis of variance propagation through the classical neural network, so they do not allow to deal with specific model under consideration. Nevertheless, we used them as a baseline strategies to compare with the proposed initialization. The proposed technique of Simple shifted initialization is described in following division. 

\subsubsection{Simple shifted initialization}

The idea behind this type of initialization is to cause minimal changes to the input signal while it is processed by the initial model. It seems a natural assumption for an untrained predistortion model, and at the same time avoids the difficulties associated with getting stuck in the domains of regularity: a model initialized in this way will always return a reasonable, although possibly not the best, result. Now, let us describe the scheme itself:
\begin{enumerate}
    \item Consider a convolution with kernel width $(2k-1)$ and weights $H = (h_1,\ldots, h_{2k-1})$. If this convolution weights are initialised as
    \begin{align*}
        h_j = 
        \begin{cases}
            1 + i,~ j = k, \\
            0,~ j\ne k.
        \end{cases}
    \end{align*}
    and zero padding of appropriate length is used, then input vector $x$ will be unchanged by applying the convolution.
    \item We initialize all entries of spline coefficient matrices with equal values $\alpha\ge 0$: $C_{ij} = \alpha,~ 1\le i, j\le P+1$.
    In our experiments, setting value $\alpha = 0.01$ showed itself as a good choice.
    \item We also use a special diagonal initialization technique for output convolutions. Namely, let model have rank $R$ and output convolution at layer $m$ have weights $H^m = (h_1^m,\ldots,h_{2k-1}^m)$. This weights are initialized as
    \begin{align*}
        h_i^j = 
        \begin{cases}
            1,~ j - i = (R+1) \text{ mod } 2 - k, \\
            0,~ \text{else}.
        \end{cases}
    \end{align*}
    Consider layer number $(R+1) \text{ mod } 2$ -- a middle layer. The corresponding output convolution at position $k$ has entry $1$ and zeros elsewhere. At the consequent layer number $(R+1) \text{ mod } 2 + 1$, the output convolution has $1$ at shifted position $k+1$. For example, if kernel width equals $2k-1=7$ and model has $R=5$ layers, the initialization is the following: $H^i = \mathbf{e}_i,\; i=1,...,5$, where $\mathbf{e}_i$ is $i$-th unit basis vector. Initializing output convolutions with ones in different positions helps make layers non-equivalent at the stage of initialization. Our experiments show that diagonal initialization significantly improves optimization effectiveness and allows to achieve better loss values.
\end{enumerate}

Numerous experiments (such as those presented in section \ref{multistart}, as well as experiments comparing the operation of other optimization methods) have shown that the presented initialization technique allows achieving better model stability in comparison with classical techniques for initializing neural networks, as well as obtaining an initial approximation with a better quality criterion value. In all the experiments presented in this article, Simple shifted initialization was used, and in no case was it observed that the problem at a given point had poor properties: the methods under consideration demonstrated satisfactory starting convergence and converged to very close solution points.

}

\section{Overfitting}
{
Overfitting is a common problem in computational graphs parameters training, when the tuned model corresponds too closely or exactly to a particular set of data, and may therefore fail to fit additional data or predict future observations reliably. In other words, in such a situation, the generalizing ability of the model is sacrificed to the quality of optimization of a specific function of empirical risk, due to which the performance of methods will decrease drastically on the data, which differ from those used for training.

In order to control overfitting, in the experiments presented in the work, the original dataset was divided into two parts~--- training (75\%) and validation (25\%), and during training, the loss function was calculated both on training and validation datasets, in order to be able to compare the quality of the methods and detect overfitting to data. Fig.~\ref{figr:overfitting} presents the results of experiments demonstrating how train and validation errors differ for several methods. As one can see, in the case of the used partition, the error difference is $0.05$ dB, whereas the error itself at the given time interval is $-37$ dB. Note that at this scale of training time and data quantity the least overfitting is achieved when using the L-BFGS method. 

\subsection{Different training set size}

From the point of view of studying the specific properties of the data generated by the signal arriving at the input of the pre-distorter, it is interesting to consider the dependence of the model's susceptibility to overfitting on the size of the training dataset. The dataset of 245~760 pairs of complex numbers $ (x, \overline {y}) $, used to set up all experiments in this work, was divided in the course of this experiment into two parts: a training set and a validation set. Moreover, the data was split in a sequential form, without random shuffling (which, on the contrary, is usually done in the case of training, for example, neural networks), so that an solid opening signal segment is used for training, and the entire remaining signal segment is taken as a validation set. After training the model on the selected training signal, the model with the resulting parameters was used to obtain a solution for the validation signal, and the performance of the method was thus assessed in parallel for two sets. The L-BFGS(900) method was used as the optimization method.

It is clearly seen from the results of the experiment that there is a discrepancy between the results of the model for different volumes of training and validation sets. With a small amount of training set, we get a model that describes very well a small number of objects (overtrained), but returns an irrelevant result when processing a new signal with naturally changed characteristics. At the same time, if you use a significant part of the dataset for training the model, the result for an arbitrary segment of the signal will, on average, be quite satisfactory, while the model can make errors for some rare individuals. Note that the specificity of this model is a rather small size of the required training sample: even when using 20\% of the training signal segment (which in this case has a size of $\approx 200~000$), the difference between the quality metrics for the training and validation samples does not exceed $0.5$ dB. At the same time, when assessing the effects of overfitting, it is important that the samples used after splitting have a sufficiently large size, since when choosing a too small training or validation set, because of data in this parts can be really different and as a result naturally occurring approximation errors begin to strongly influence the result. 

\vspace{-0.8cm}
\begin{minipage}[t]{\dimexpr.4\textwidth-.5\columnsep}
\raggedright

\begin{figure}[H]
\centering
\includegraphics[width=\linewidth]{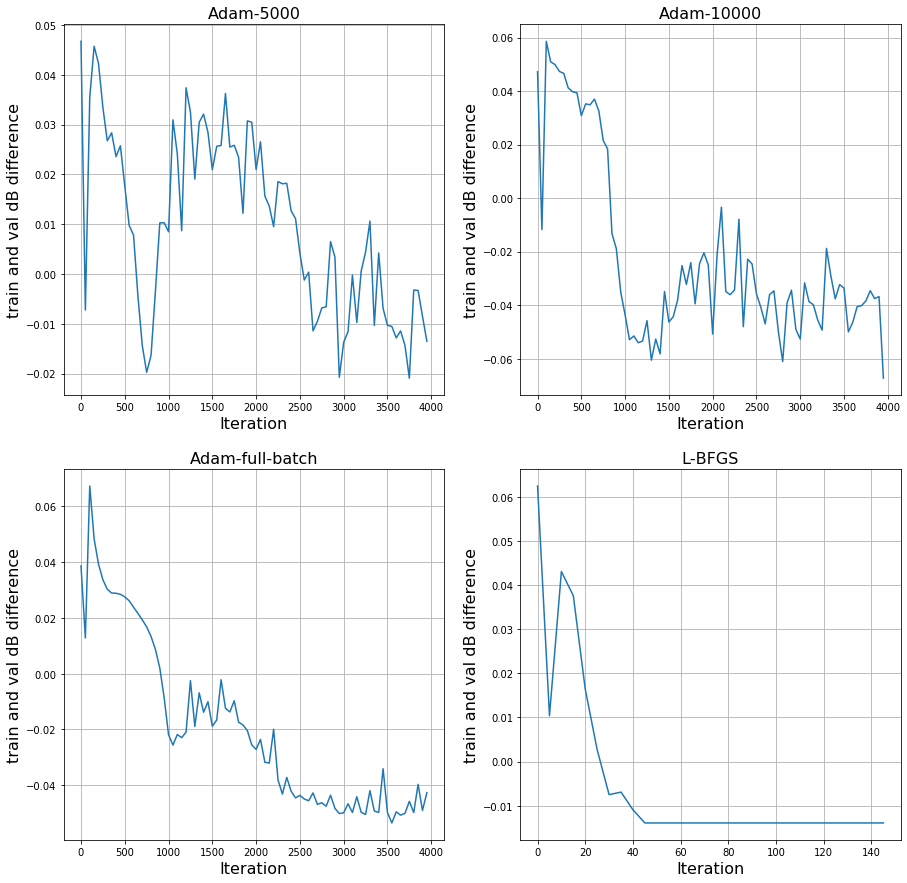}
\caption{Difference between train and validation errors}
\label{figr:overfitting}
\end{figure}

\end{minipage}
\begin{minipage}[t]{\dimexpr.55\textwidth-.5\columnsep}
\raggedleft

\begin{figure}[H]
\centering
\includegraphics[width=0.35\linewidth]{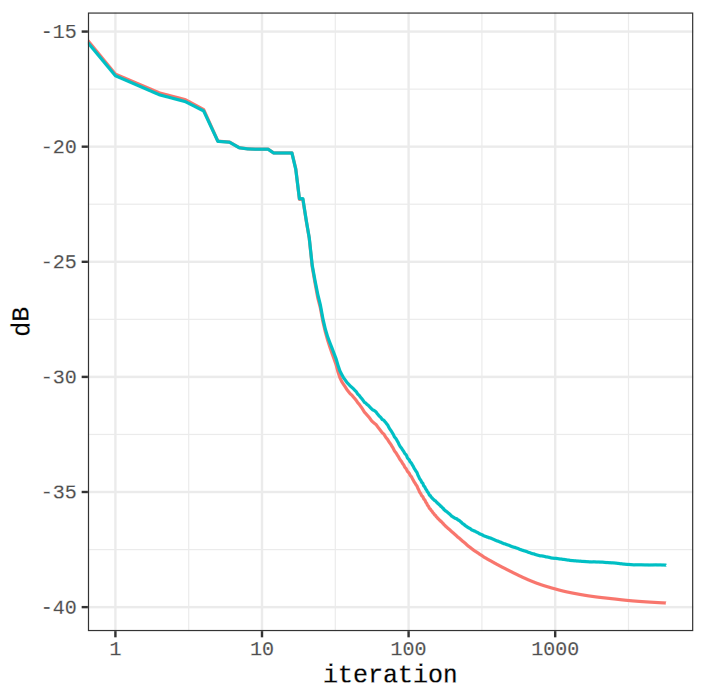}
\includegraphics[width=0.48\linewidth]{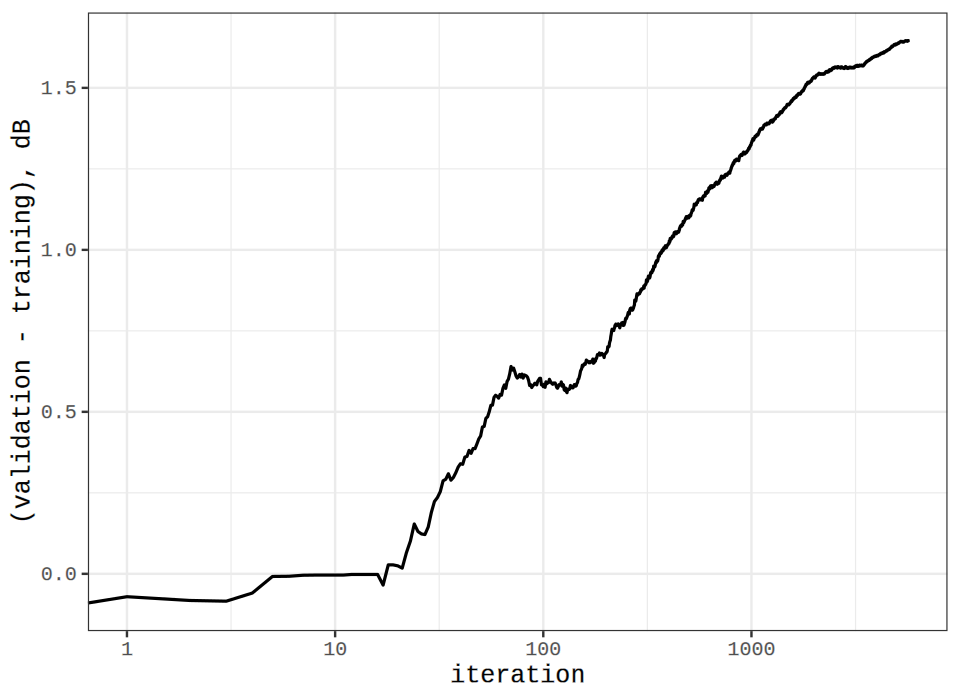}
\hspace{3cm}
\includegraphics[width=0.35\linewidth]{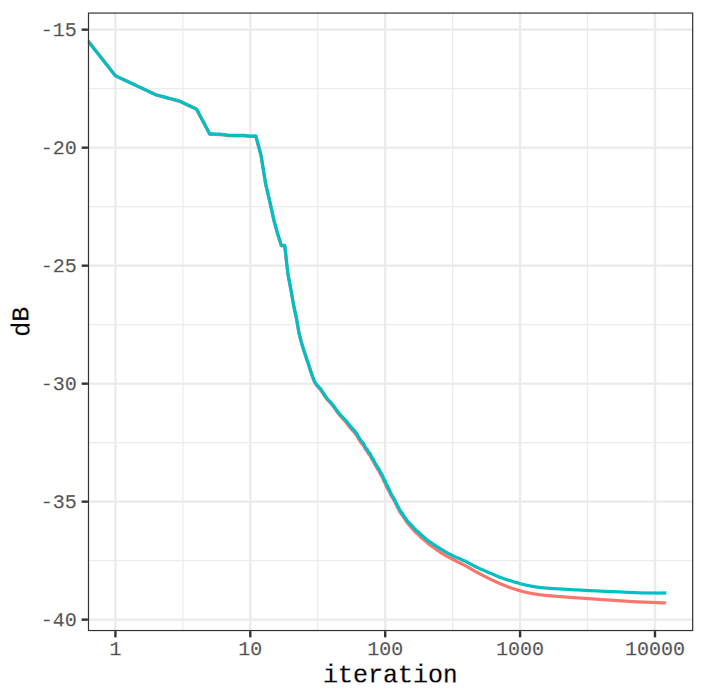}
\includegraphics[width=0.48\linewidth]{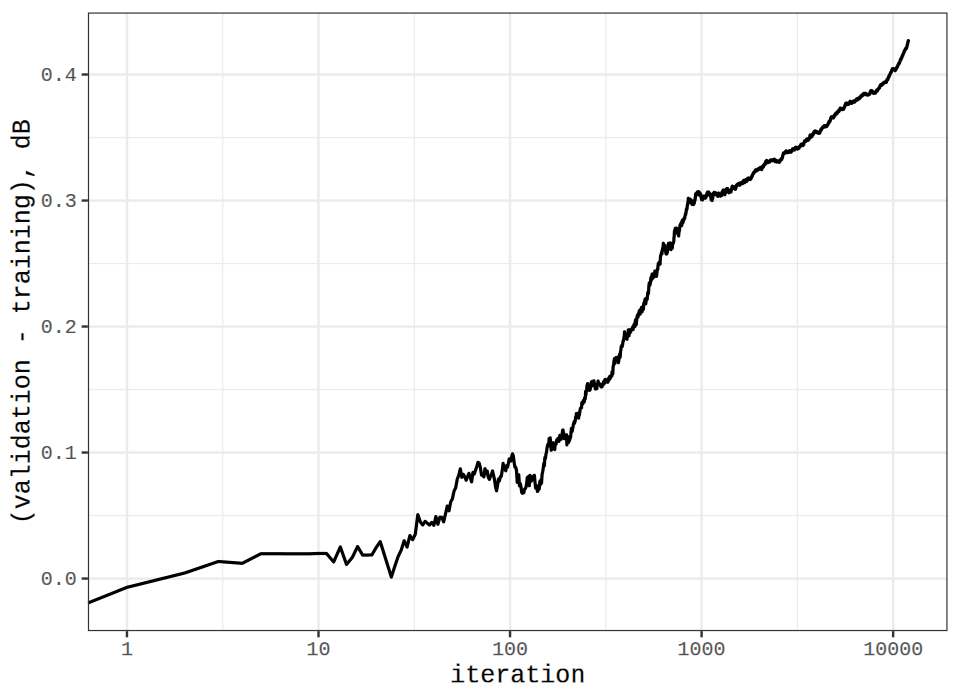}
\caption{Residual model overfitting, 5\%, 20\% of original data used as training set. Red line~--- convergence while training on sample signal, blue line~--- quality of the solution for validation signal}
\label{fig:overfitting_99_1}
\end{figure}

\end{minipage}
\vspace{-0.6cm}
}

\section{Conclusion}
{

This article discusses various approaches to optimizing the parameters of computational graphs simulating the behavior of a digital pre-distorter for the modulated signal. 
In the numerous experiments it was tested different full-gradient methods, and stochastic algorithms.

Among the many randomized (batch) algorithms that significantly use the sum-type structure of the objective functional, the Adam algorithm, which is most relevant for use in online-training regime, demonstrates the best efficiency. However, it should be noted that for its effective adjustment, that is, choosing the optimal step length and batch size, it is necessary to perform a rather time-consuming pre-calculations.

Of all the considered methods, the L-BFGS algorithm turned out to be the undisputed leader. It may be somewhat unexpected that the optimal memory depth of the L-BFGS method for this problem is in the range of 800--1000. Note that the idea of using the L-BFGS method for DPD optimization was proposed earlier in paper \cite{bao2013restarted}, for models based on the Volterra series. Experiments described in this paper thus confirm the particular effectiveness of the L-BFGS method for the DPD problem in relative independence from a particular model and dataset. At the same time, the idea of deep memory is original, and perhaps it is specific to the used model class. In addition to the best convergence rate, the L-BFGS method as a training procedure also leads to the least error on the validation set: the discrepancy in the quality metric, that characterizes the overfitting susceptibility, for the dataset used is approximately $0.05$ dB.

There were described a number of new modifications of the Gauss--Newton method proposed by Yu.E. Nesterov, including the Method of Stochastic Squares. The practical efficiency of the proposed approaches is not only significantly higher than that of other Gauss--Newton methods, but is also the best among all the local methods considered in the experiments (see fig. \ref{fig:ssm_conv}).

Many experiments have been carried out evaluating the specifics of the dataset used, generated by the samples of the modulated signal. Experiments on the use of different sizes of the training sample have shown that it is enough to use 20\% of the original dataset to obtain a sufficiently good quality on the validation set ($-38$ dB). Moreover, in this case, it is possible to reach the $-38$ dB threshold much faster than using the full training dataset. It should be noted that even 5\% of the data is enough to reach the $-37$ dB threshold, and also in a much shorter time.


\vspace{-0.6cm}
\hspace{-1cm}
\begin{minipage}[t]{\dimexpr.6\textwidth-.5\columnsep}
\raggedright
\begin{table}[H]
\caption{Best methods performance, residual model}
\centering
\begin{tabular}{l||r|r|r|r}
\multirow{2}{*}{Method} & \multicolumn{4}{c}{Time to reach dB level, sec.} \\ \cline{2-5}
& -30 dB & -35 dB & -37 dB & -39 dB \\ \hline \hline

SDM & 25.83 & 925.92 & 5604.51 & \\

CG(DY) & 12.38 & 47.21 & 80.74 & 3721.43 \\

DFP(inf) & 11.10 & 43.96 & 72.56 & 944.52 \\
BFGS(100) & \ok{30}{7.41} & 34.77 & 53.89 & 695.28 \\

LBFGS(900) & \ok{70}{4.01} & 16.75 & 34.49 & \ok{40}{399.52} \\

3SM & 633.88 & 1586.69 & 1747.28 & 4616.80 \\
SSM($6 n$) & 47.32 & 72.30 & 89.08 & \ok{60}{314.95} \\

\end{tabular}
\vspace{0.5cm}
\end{table}

\end{minipage}
\begin{minipage}[t]{\dimexpr.5\textwidth-.5\columnsep}
\raggedleft

\begin{figure}[H]
\centering
\includegraphics[width=\linewidth]{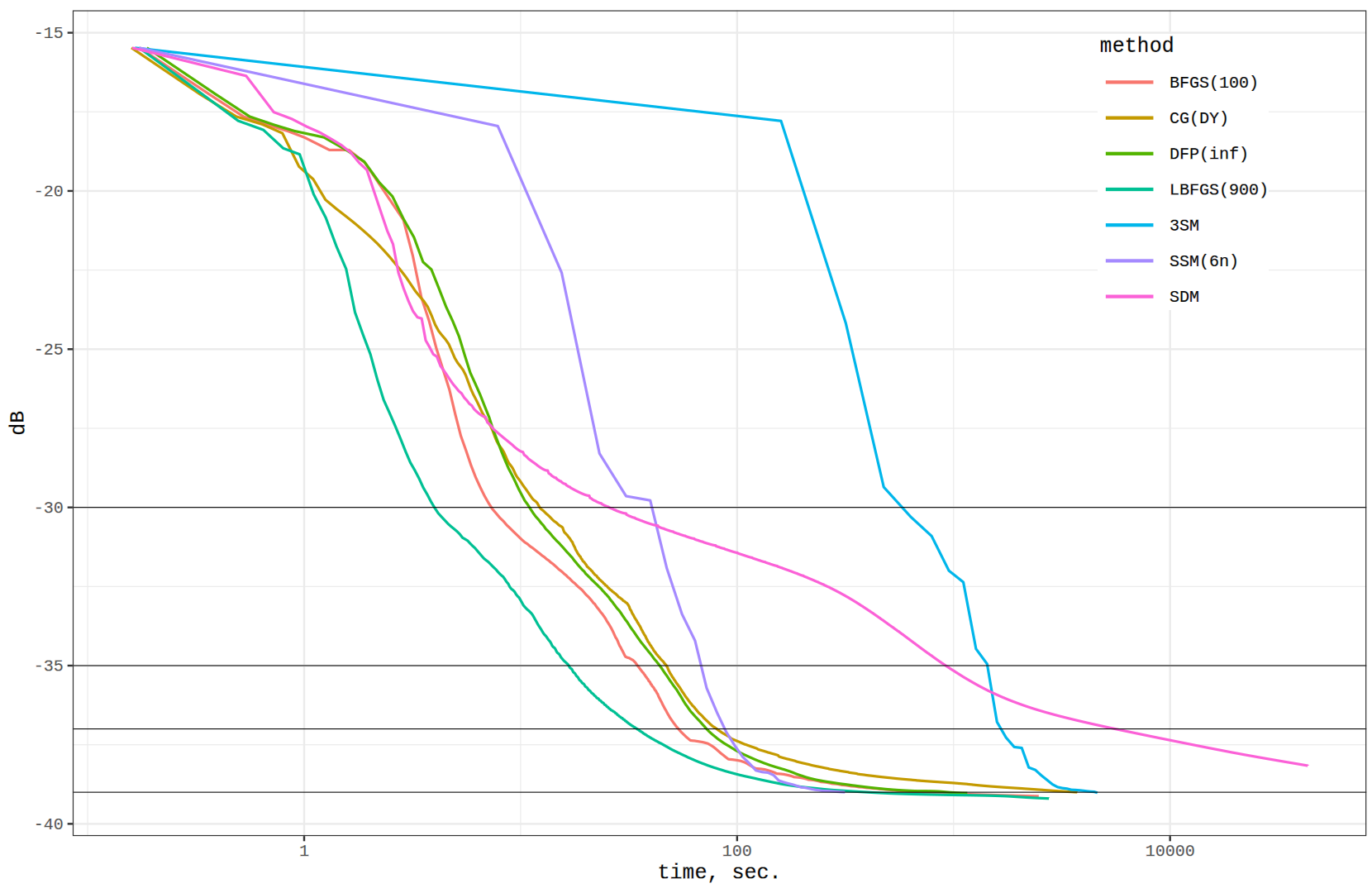}
\caption{Best methods convergence, residual model}
\end{figure}

\end{minipage}
\vspace{-0.1cm}

Thus, despite the fact that the considered class of models has significant specificity, following the classical way of studying large computational graphs from the point of view of their parameters optimizing, formed mainly around the problem of training neural networks, makes it possible to collect a set of algorithms and approaches that are most effective for the problem under consideration. Moreover, many of the solutions developed specifically for neural networks turn out to be relevant for Wiener--Hammerstein type models. In particular, adaptive stochastic methods remain just as effective. At the same time, it is possible to significantly and drastically improve the results of classical approaches taking into account the specifics of the problem, such as the use of deep memory for the L-BFGS method, shifted weights initialization, small width of network layers or a small training sample. Apparently, the dependencies found in the course of the described study are quite universal for this family of models, and therefore the presented observations can be useful not only for efficiently solving related practical problems, but also for further exploration of the problem of digital predistortion of the signals. 
}

\bibliographystyle{splncs04.bst}
\bibliography{bibl}

\begin{thebibliography}{10}
\providecommand{\url}[1]{\texttt{#1}}
\providecommand{\urlprefix}{URL }
\providecommand{\doi}[1]{https://doi.org/#1}

\bibitem{amir2021sgd}
Amir, I., Koren, T., Livni, R.: Sgd generalizes better than gd (and
  regularization doesn't help). arXiv preprint arXiv:2102.01117  (2021)

\bibitem{auer1996exponentially}
Auer, P., Herbster, M., Warmuth, M.K.: Exponentially many local minima for
  single neurons. In: Advances in neural information processing systems. pp.
  316--322 (1996)

\bibitem{bao2013restarted}
Bao, J.L., Zhu, R.X., Yuan, H.X.: Restarted lbfgs algorithm for power amplifier
  predistortion. In: Applied Mechanics and Materials. vol.~336, pp. 1871--1876.
  Trans Tech Publ (2013)

\bibitem{barzilai1988}
Barzilai, J., Borwein, J.M.: Two-point step size gradient methods. IMA Journal
  of Numerical Analysis  \textbf{8},  141--148 (1988)

\bibitem{choromanska2015loss}
Choromanska, A., Henaff, M., Mathieu, M., Arous, G.B., LeCun, Y.: The loss
  surfaces of multilayer networks. In: Artificial intelligence and statistics.
  pp. 192--204 (2015)

\bibitem{dennis1977quasi}
Dennis, Jr, J.E., Mor{\'e}, J.J.: Quasi-newton methods, motivation and theory.
  SIAM review  \textbf{19}(1),  46--89 (1977)

\bibitem{garipov2018loss}
Garipov, T., Izmailov, P., Podoprikhin, D., Vetrov, D., Wilson, A.G.: Loss
  surfaces, mode connectivity, and fast ensembling of dnns. arXiv preprint
  arXiv:1802.10026  (2018)

\bibitem{ghannouchi2015behavioral}
Ghannouchi, F.M., Hammi, O., Helaoui, M.: Behavioral modeling and predistortion
  of wideband wireless transmitters. John Wiley \& Sons (2015)

\bibitem{griewank1989automatic}
Griewank, A., et~al.: On automatic differentiation. Mathematical Programming:
  recent developments and applications  \textbf{6}(6),  83--107 (1989)

\bibitem{hajek1988cooling}
Hajek, B.: Cooling schedules for optimal annealing. Mathematics of operations
  research  \textbf{13}(2),  311--329 (1988)

\bibitem{haykin2008adaptive}
Haykin, S.S.: Adaptive filter theory. Pearson Education India (2008)

\bibitem{he2016deep}
He, K., Zhang, X., Ren, S., Sun, J.: Deep residual learning for image
  recognition. In: The IEEE Conference on Computer Vision and Pattern
  Recognition (CVPR) (06 2016)

\bibitem{kirkpatrick1983optimization}
Kirkpatrick, S., Gelatt, C.D., Vecchi, M.P.: Optimization by simulated
  annealing. science  \textbf{220}(4598),  671--680 (1983)

\bibitem{liu1989limited}
Liu, D.C., Nocedal, J.: On the limited memory bfgs method for large scale
  optimization. Mathematical programming  \textbf{45}(1),  503--528 (1989)

\bibitem{marquardt1963algorithm}
Marquardt, D.W.: An algorithm for least-squares estimation of nonlinear
  parameters. Journal of the society for Industrial and Applied Mathematics
  \textbf{11}(2),  431--441 (1963)

\bibitem{andrei40cg2008}
Neculai, A.: Conjugate gradient algorithms for unconstrained optimization. a
  survey on their definition. ICI Technical Report  \textbf{13},  1--13 (2008)

\bibitem{nesterov2021flexible}
Nesterov, Y.: Flexible modification of gauss-newton method. CORE Discussion
  paper  (2021)

\bibitem{nocedal1980updating}
Nocedal, J.: Updating quasi-newton matrices with limited storage. Mathematics
  of computation  \textbf{35}(151),  773--782 (1980)

\bibitem{nocedal2006numerical}
Nocedal, J., Wright, S.: Numerical optimization. Springer Science \& Business
  Media (2006)

\bibitem{nolte2000note}
Nolte, A., Schrader, R.: A note on the finite time behavior of simulated
  annealing. Mathematics of Operations Research  \textbf{25}(3),  476--484
  (2000)

\bibitem{polyak1969}
Polyak, B.T.: Minimization of unsmooth functionals. USSR Computational
  Mathematics and Mathematical Physics  \textbf{9},  14--29 (1969)

\bibitem{price2006differential}
Price, K., Storn, R.M., Lampinen, J.A.: Differential evolution: a practical
  approach to global optimization. Springer Science \& Business Media (2006)

\bibitem{schreurs2008rf}
Schreurs, D., O'Droma, M., Goacher, A.A., Gadringer, M.: RF power amplifier
  behavioral modeling. Cambridge university press New York, NY, USA (2008)

\bibitem{skajaa2010limited}
Skajaa, A.: Limited memory bfgs for nonsmooth optimization. Master's thesis
  (2010)

\bibitem{storn1995differrential}
Storn, R.: Differrential evolution-a simple and efficient adaptive scheme for
  global optimization over continuous spaces. Technical report, International
  Computer Science Institute  \textbf{11} (1995)

\bibitem{yudin2021flexible}
Yudin, N., Gasnikov, A.: Flexible modification of gauss-newton method and its
  stochastic extension. arXiv preprint arXiv:2102.00810  (2021)

\bibitem{zhang2019adam}
Zhang, J., Karimireddy, S.P., Veit, A., Kim, S., Reddi, S.J., Kumar, S., Sra,
  S.: Why adam beats sgd for attention models. arXiv preprint arXiv:1912.03194
  (2019)

\end{thebibliography}

\vspace{5.6cm}

\appendix
\addcontentsline{toc}{section}{Appendices}
\section*{Appendices}
\vspace{1.6cm}

\section{Stochastic Squares Method Implementation Details}
The efficiency of Stochastic Squares Method implementation heavily depends on the size $p$ of the batch. We follow the work \cite{nesterov2021flexible}.

\subsection{Small batch size (\texorpdfstring{$0 \leq p \leq n$}{TEXT})} 

At each iteration of this method, it is necessary to invert the
matrix
\[
 B_k := {\frac{1}{2 p \hat f_1(x)}} G_k G_k^T + 2^{i_k} L_k I \in \R^{n \times n}.
\]
If $p$ is small, then it is acceptable to do this at each iteration using Sherman--Morrison formula:
\[
 B_k^{-1} = {\frac{1}{2^{i_k} L_k}} \left[I + \gamma_k G_k G_k^T\right]^{-1} \; = \; {\frac{1}{2^{i_k} L_k}} \left[I - \gamma_k G_k [I_p + \gamma_k G_k^T G_k ]^{-1} G_k^T \right],
\]
where $\gamma_k = {\frac{1}{2^{i_k+1} p L_k \hat f_1(x)}}$. Note that we need $B_k^{-1}$ only for multiplying by vectors. Therefore, we have to form and invert only a small matrix
\[
 [I_p + \gamma_k G_k^T G_k] \in \R^{p \times p}.
\]
This requires $\mathcal{O}(p^2 (n+p))$ operations, and this must be done at all $i_k$ internal iterations of step d) in method. Note that the version with $p=0$ corresponds to the simple Gradient Method as applied to function $\hat f_1$.

\subsection{Big batch size (\texorpdfstring{$n < p \leq m$}{TEXT})} 

In this case, $n$ is the smallest size of the matrix $G_k$. Therefore, for efficient implementation of the inversion of matrix $B_k$, it is reasonable to compute an appropriate factorization of matrix $G_k \in \R^{n \times p}$. Namely, we need to represent this matrix in the following form:
\[
 G_k = U_k \Lambda_k V_k,
\]
where $U_k$ and $V_k$ are orthogonal matrices of the corresponding sizes and $\Lambda_k \in \R^{n \times n}$ is the lower-triangular bi-diagonal matrix. It is easier to compute this factorization if we represent the matrices $U_k$ and $V_k$ in multiplicative form. The basic orthogonal transformation used in this process is the Hausdorff matrix
\[
 H(a,b) = I - 2 {\frac{(b-a)(b-a)^T}{\| b - a \|_2^2}},
\]
where $b \neq a$ and $\| a \|_2 = \| b \|_2$. It is easy to check that
\[
 H(a,b)a = b, \quad H(a,b)b \; = \; a, \quad H^2(a,b) \; = \; I.
\]

The orthogonalization process has $n$ steps. At the first step, we multiply the matrix $G_k$ from the right by an orthogonal matrix $S_1 \in \R^{p \times p}$, which maps the first row $e_1^\top G_k$ to the vector $\| G_k^\top e_1 \|_2\; e_1^\top$, where $e_i$ is the $i$-th unit basis vector of the dimension corresponding to the context.

Further, denote by $G_k^+$ the matrix composed by the rows $2 \dots n$ of the matrix $G_k S_1$. We multiply it from the left by an orthogonal matrix, which transforms its first column into a vector proportional to $e_1$. This ends the first step of the process. Denote the resulting matrix by $G_k^{++}$. The matrix for the next step is formed by columns $2 \dots n$ of matrix $G_k^{++}$.

In the end, we get a lower-triangular bi-diagonal matrix $\Lambda_k$ and the matrices $U_k$ and $V_k$ represented in multiplicative form. For the optimization process, it is reasonable to keep them in this form since we need only to multiply these matrices by vectors.

For stability of the above process, it is important to apply matrices $H(a,b)$ with sufficiently big $\| a - b \|_2$. This can be always achieved by considering vectors $\pm b$ (for us they are multiples of the coordinate vectors), and choosing the sign ensuring $\la a, b \ra \leq 0$.

\section{Quasi-Newton Algorithms}
\vspace{-1cm}
\begin{table}[H]
\centering
\caption{BFGS/DFP algorithmic parameters}
\begin{tabularx}{\linewidth}{l|X}
  $\varepsilon_{GN}$ & the accuracy of the stopping criterion according to the gradient norm \\ \hline

  $K_{RS}$ & the direction restart frequency \\
\end{tabularx}
\end{table}
\vspace{-1cm}

\begin{algorithm}[H]
	\caption{BFGS/DFP method}
  \label{alg:BFGS_DFP}
	\begin{algorithmic}[1]
		\REQUIRE algorithmic parameters setup, $x^{0}$.
		\FOR{$k = 0, \hdots, k_{\max}$}
		
		\IF	{$ \| \nabla f(x^{k}) \| \le \varepsilon_{GN}$}
  		\STATE $x^* = x^{k}$
  		\STATE \textbf{return}
		\ENDIF
		
		\IF {($k \mod K_{RS}) = 0$}
  		\STATE $H^k = E$
		\ENDIF
		
		\STATE $p^{k} = -\nabla f(x^{k}) H^{k}$
		\STATE $\alpha^{k} = \argmin_{\alpha > 0} f(x^{k} + \alpha p^{k})$.
		\STATE $x^{k+1} = x^{k} +\alpha^{k} p^{k}$
		\STATE $H^{k+1} = \Phi(Z^{k}, x^{k}, x^{k+1}, \nabla f(x^{k}), \nabla f(x^{k+1}))$ -- update $H$ ($B^{-1}$ approximation) with corresponding update rule $\Phi$ for BFGS/DFP scheme
		\ENDFOR
		\ENSURE $x^*$ if the early termination condition was reached, otherwise the final iterate $x^{(k_{\max} + 1)}$.
	\end{algorithmic}
\end{algorithm}


\vspace{5.2cm}
\begin{table}[h]
\centering
\caption{L-BFGS algorithmic parameters}
\begin{tabularx}{\linewidth}{l|X}
  $m > 0$ & the history size (length). \\ \hline

  $\varepsilon_{GN}$ & the accuracy of the stopping criterion according to the gradient norm \\ \hline

  $\varepsilon_{FLT}$ & ``machine epsilon'', equal to $7 \cdot 10^{-16}$ ($1 + \varepsilon_{MAS} = 1$).\\ \hline
\end{tabularx}
\end{table}

\vspace{-1.2cm}
\begin{algorithm}[H]
	\caption{L-BFGS method}
  \label{alg:lbfgs}
	\begin{algorithmic}[1]
		
		\REQUIRE algorithmic parameters setup, initial $x^{0}$.
		
    \STATE $p^{0} =-\nabla f(x^{0})$.
		\FOR{$k = 0, \hdots, k_{\max}$}
		\IF	{$ \| \nabla f(x^{k}) \|_2 \le \varepsilon_{GN}$}
  		\STATE $x^* =  x^{k}$
  		\STATE \textbf{return}
		\ENDIF		
		
    \STATE $\alpha^k = \argmin_{\alpha > 0} f(x^k + \alpha p^k)$
    \STATE $x^{k+1} = x^{k} +\alpha^{k} p^{k}$
		\STATE $p^{k+1} = -\nabla f(x^{k+1})$.

		\STATE $\tilde{m} = \min(k + 1, m)$
    \FOR {$j=\tilde{m}, \tilde{m}-1, \hdots, 2$}
      \STATE $s^{j} = s^{j-1}$.
      \STATE $y^{j} = y^{j-1}$.
    \ENDFOR
		\STATE $s^{1} = x^{k+1} - x^{k}$
    \STATE $y^{1} = \nabla f(x^{k+1} )-\nabla f(x^{k})$

		\FOR {$j=1,2, \hdots, \tilde{m}$}
		\STATE  $\rho = \left\langle s^{j}, y^{j} \right\rangle$
		\IF {$ \left|\rho \right| \le \varepsilon_{FLT}$}
    \STATE  $p^{k + 1} = -\nabla f(x^{k + 1})$, \textbf{go to line~\ref{lbfgs:st2}}
    \ENDIF
    \STATE $\hat{\alpha}^{j} = 1 / \rho \left\langle s^{j} ,p^{k + 1} \right\rangle $.
    \STATE $p^{k + 1} = p^{k + 1} - \hat{\alpha}^{j} y^{j}$.
		\ENDFOR
		
		\FOR {$j=\tilde{m},\tilde{m}-1, \hdots, 1$}
		\STATE $\rho = \left\langle s^{j}, y^{j} \right\rangle$
		\IF {$ \left|\rho \right| \le \varepsilon_{FLT}$}
      \STATE  $p^{k+1} =-\nabla f(x^{k+1} )$, \textbf{go to line~\ref{lbfgs:st2}}
    \ENDIF
    \STATE $\hat{\beta}^{j} = 1 / \rho \left\langle y^{j} ,p^{k+1} \right\rangle $. 
    \STATE $p^{k+1} = p^k + (\hat{\alpha}^{j} - \hat{\beta}^{j})s^{j} $.
		\ENDFOR
		\STATE \textbf{continue loop} \label{lbfgs:st2}
		\ENDFOR
		\ENSURE $x^*$ if the early termination condition was reached, otherwise the final iterate $x^{(k_{\max} + 1)}$.
		
	\end{algorithmic}
\end{algorithm}
\section{Conjugate Gradients Algorithms}
\vspace{-1cm}
\vspace{.1cm}
\begin{table}[H]
\centering
\caption{CG algorithmic parameters}
\begin{tabularx}{\linewidth}{l|X}
  $\varepsilon_{GN}$ & the accuracy of the stopping criterion according to the gradient norm \\ \hline

  $\varepsilon_{FLT}$ & ``machine epsilon'', equal to $7 \cdot 10^{-16}$ ($1 + \varepsilon_{MAS} = 1$).\\ \hline

  $K_{RS}$ & the direction restart frequency \\
\end{tabularx}
\end{table}
\vspace{-1.4cm}

\begin{algorithm}[H]
	\caption{CG(HS) method}\label{alg:CG_HS}
	\begin{algorithmic}[1]
		\REQUIRE algorithmic parameters setup, initial $x^{0}$.
    \STATE $p^{0} =-\nabla f(x^{0})$
		\FOR{$k = 0, \hdots, k_{\max}$}
		\IF	{$ \| \nabla f(x^{k}) \| \le \varepsilon_{GN}$}
		  \STATE $x^* = x^{k}$
		  \STATE \textbf{return}
		\ENDIF

    \STATE $\alpha^k = \argmin_{\alpha > 0} f(x^k + \alpha p^k)$
    \STATE $x^{k+1} = x^{k} +\alpha^{k} p^{k}$

    \IF {($(k + 1) \mod K_{RS}) = 0$}
      \STATE  $p^{k+1} =-\nabla f(x^{k+1} )$, \textbf{go to line~\ref{cg:end_loop}}
		\ENDIF

		\STATE $s = \langle \nabla f(x^{k+1}) - \nabla f(x^{k}), x^{k+1} - x^{k} \rangle$
		\IF {$\left| s \right| \le \varepsilon_{FLT}$}
      \STATE  $p^{k+1} =-\nabla f(x^{k+1} )$, \textbf{go to line~\ref{cg:end_loop}}
    \ENDIF

    \STATE $\beta = \langle \nabla f(x^{k+1}) - \nabla f(x^{k}), \nabla f(x^{k}) \rangle / s$
		
		\STATE $p^{k+1} =-\nabla f(x^{k+1} ) + \beta \cdot p^{k} $
    \STATE \textbf{continue loop} \label{cg:end_loop}
		\ENDFOR
    \ENSURE $x^*$ if the early termination condition was reached, otherwise the final iterate $x^{(k_{\max} + 1)}$.
	\end{algorithmic}
\end{algorithm}
\vspace{-1.5cm}
\begin{algorithm}[H]
	\caption{CG(FR) method}\label{alg:CG_FR}
	\begin{algorithmic}[1]
		\REQUIRE algorithmic parameters setup, initial $x^{0}$.
    \STATE $p^{0} =-\nabla f(x^{0})$
		\FOR{$k = 0, \hdots, k_{\max}$}
		\IF	{$ \| \nabla f(x^{k}) \| \le \varepsilon_{GN}$}
		  \STATE $x^* = x^{k}$
		  \STATE \textbf{return}
		\ENDIF

    \STATE $\alpha^k = \argmin_{\alpha > 0} f(x^k + \alpha p^k)$
    \STATE $x^{k+1} = x^{k} +\alpha^{k} p^{k}$

    \IF {($(k + 1) \mod K_{RS}) = 0$}
      \STATE  $p^{k+1} =-\nabla f(x^{k+1} )$
    \ELSE
      \STATE $\beta= || \nabla f(x^{k+1} )||^2 \  / \  || \nabla f(x^{k} )||^2$
      \STATE $p^{k+1} =-\nabla f(x^{k+1} ) + \beta \cdot p^{k} $
    \ENDIF
    \ENDFOR
    \ENSURE $x^*$ if the early termination condition was reached, otherwise the final iterate $x^{(k_{\max} + 1)}$.
	\end{algorithmic}
\end{algorithm}

\begin{algorithm}[H]
\caption{CG(Nesterov) method}
  \label{alg:CG_Nest}
	\begin{algorithmic}[1]
		\REQUIRE algorithmic parameters setup, initial $x^{0}$.
    $y^{-1} = y^{-2} = x^0$
		\FOR{$k = 0, \hdots, k_{\max}$}
		  \IF {$ \| \nabla f(x^{k}) \| \le \varepsilon_{GN}$}
		    \STATE $x^* = x^{k}$
		    \STATE \textbf{return}
		  \ENDIF

      \STATE $\alpha_1^{k} =\argmin_{\alpha > 0} f(x^{k} + \alpha (y^{k-2} - x^{k}))$
      \STATE $y^{k} = x^{k} + \alpha_1^k (y^{k-2} - x^{k})$

		  \STATE $\alpha_2^{k} =\argmin_{\alpha > 0} f(y^{k} - \alpha \nabla f(y^k))$
		  \STATE $x^{k+1} = y^{k} - \alpha_2^k \nabla f(y^k)$
		\ENDFOR
    \ENSURE $x^*$ if the early termination condition was reached, otherwise the final iterate $x^{(k_{\max} + 1)}$.
	\end{algorithmic}
\end{algorithm}
\vspace{-1.2cm}

\vspace{10cm}
\section{Global Algorithms}\label{global}

\subsection{Simulated annealing}

\begin{table}[H]
\centering
\vspace{-1cm}
\caption{Simulated annealing algorithmic parameters}
\label{sa_params}

\begin{tabularx}{\linewidth}{l|X}
  ${t_0}$ & initial "temperature" $[1.0, 1000.0]$, $6.0$ \\ \hline

  $K_{CG}$ & number of iterations of the conjugate gradient method $[1, 1000]$ $0$ (v1) or $50$ (v2) \\ \hline

  $K_{jump}$ & frequency of "jump" -- forced update of the starting point $[1, 1000]$, $100$ \\ \hline

  $D_{jump}$ & "jump" intensity $[10^{-6}, 10^3]$, $1$ \\ \hline

  $t_{\max}$ & algorithm time limit, sec. $[1, 10^5]$, $600$ \\
\end{tabularx}
\end{table}

\begin{algorithm}[H]
  \caption{Simulated Annealing method}
  \label{alg:sa}
  \begin{algorithmic}[1]
    \REQUIRE algorithmic parameters setup, initial $x^0$, bounds ${x_L}, {x_G} \in \mathbb{R}^n$.

    \STATE $x^{*} = x^0$
    \STATE $f^{*} = f(x^{*})$
    \STATE $k = 0$

    \WHILE{$t^k \le t_{\max}$}
      \STATE Perform $K_{CG}$ iterations of the conjugate gradient method which begins in $x^{k}$, save results to $\bar x^k$
      \IF {$f(\bar x^k) < f(x^k)$}
        \STATE $x^{*} = \bar x^k$
        \STATE $x^{k + 1} = \bar x^k$
      \ELSE
        \STATE Generate (uniform) random vector $r \in \mathbb{R}^n: R_i \in [-1, 1]$ \label{sa:jump}
        \STATE Perform a ``jump'': $x^{k + 1} = x^{*} + D_{jump} \cdot r$
        \STATE $k = k + 1$, \textbf{goto next iteration}
      \ENDIF
      
      \STATE Calculate the ``jump'' probability $P^k = (1/k)^{1/{t_0}}$.
      \STATE Generate (uniform) random value $p \in [0, 1]$
      \IF {($K \mod K_{jump} = 0)$ \textbf{or} $(p \le P^K)$}
        \STATE Perform a ``jump'' (see line~\ref{sa:jump} for details)
      \ELSE
        \STATE $x^{k+1} = x^k$
      \ENDIF
      
      \STATE $k = k + 1$
    \ENDWHILE
  \RETURN $x^{*}$
	\end{algorithmic}
\end{algorithm}
\vspace{-1.0cm}

\subsection{Differential Evolution results}
\vspace{-1cm}
\begin{table}[H]
\centering
\caption{Differential evolution algorithmic parameters}
\begin{tabularx}{\linewidth}{l|X}
  $m$ & number of individuals in a population $[4, 1000]$, $20$ \\ \hline

  $K_{CG}$ & number of iterations of the conjugate gradient method $[1, 1000]$ $0$ (v1) or $50$ (v2) \\ \hline

  $F$ & mutation ``force'' -- amplitude of the introduced disturbance $[1.0, 1000.0]$, $0.5$ \\ \hline

  $CR$ & mutation probability $[10^{-3}, 1.0]$, $0.1$ \\ \hline
  
  $\varepsilon_{BIO}$ & minimum allowed measure of ``biodiversity'' $[10^{-12}, 10.0]$, $10^{-6}$ \\ \hline

  $t_{\max}$ & algorithm time limit, sec. $[1, 10^5]$, $600$ \\
\end{tabularx}
\end{table}

\begin{algorithm}[H]
  \caption{Differential Evolution method}
  \label{alg:diff_evo}
  \begin{algorithmic}[1]
    \REQUIRE algorithmic parameters setup, initial $x^0$, bounds ${x_L}, {x_G} \in \mathbb{R}^n$.

    \STATE $x^{*} = x^0$
    \STATE $f^{*} = f(x^{*})$
    \STATE $k = 0$

    \STATE Generate a random starting population ${P^0} = \{ {p_1} = {x^0},{p_2},...,{p_m}\}$, $x_L \le p_i \le x_G, \forall i \in [1,...,m]$

    \STATE Perform a selection of the starting population: make $m$ local descents, each of $K_{CG}$ iterations of the conjugate gradient method, each begins with its individual of the starting population $P^0$, save results to $S = \{ {s_1},{s_2},...,{s_m}\}$
    
    \STATE Refine the record: $x^{*} = \argmin\limits_{x \in S} f(x)$
    
    \STATE $P^0 = S$

    \WHILE{$t^k \le t_{\max}$}
      \STATE Calculate a weak "measure of biodiversity" $M_{BIO} = 1/(m - 1) \sum\limits_{l = 1}^{m - 1} {\| p_{l + 1} - p_l \|}$
      
      \IF{$M_{BIO} \le \varepsilon_{BIO}$}
        \RETURN $x^{*}$
      \ENDIF
      
      \FOR{$l = 1, \hdots, m$}
        \STATE Generate (uniform) random value ${R_l} \in [0,1]$
        \IF {${R_l} < CR$}
          \STATE $p_l^{k + 1} = p_l^k$
        \ELSE
          \STATE Generate three random indexes $j_1, j_2, j_3 : \  j_1 \ne j_2 \ne j_3 \ne l$
          \STATE Generate an individual of a new population 
          $p_l^{k + 1} = p_{j_1} + F (p_{j_2} - p_{j_3})$.
        \ENDIF
        \STATE Estimate new individual: perform $K_{CG}$ iterations of the conjugate gradient method which begins in $p_l^{k + 1}$, save results to ${s_l}$
      \ENDFOR

      \STATE Refine the record: $x^{*} = \argmin\limits_{x \in S} f(x)$

      \STATE Perform selection: add $m$ best individuals from $P^k \cup S$ to the new population ${P^{k + 1}}$
      
      \STATE $k = k + 1$
    \ENDWHILE
  \RETURN $x^{*}$
	\end{algorithmic}
\end{algorithm}
\vspace{-1.2cm}


\section{Other Algorithms}
\subsection{Raider Method}

The main idea of the method is to identify a subset of influential variables and to optimize at each iteration in the gradient direction, corresponding only to this set. To assess the influence of variables, an analysis of the gradient component modules is used, from the maximum of which a cross section is performed, which determines the influence of each variable in this situation. The cross-sectional level is set by the algorithmic parameter $D_{LEVEL}$. In the direction thus truncated, the simplest one-dimensional search is performed, which implements a multiple division of the step (starting from a single value) to achieve an improving approximation. It seems obvious that the effectiveness of such a method should most significantly depend on the value of the parameter , indicating which part of the variables is selected at a given iteration. In connection with the above, we consider two versions of the algorithm: the ``good'' variant with  $D_{LEVEL} =0.2$ (Var1) and the ``evil'' variant with $D_{LEVEL} = 0.9$ (Var2).

\begin{table}[H]
\centering
\caption{Raider algorithmic parameters}
\begin{tabularx}{\linewidth}{l|X}
  $D_{LEVEL}$ & the gradient cross section level $[10^{-6}, 1.0]$, equal to $0.2$. \\ \hline

  $\varepsilon_{NG}$ & the accuracy of the stopping criterion according to the gradient norm from $[10^{-12}, 10^2]$, equal to $10^{-5}$.\\ \hline



  $K_{\alpha}$ & the coefficient of the step compression from $[1.1, 10.0]$, equal to $3.0$ (Var2).\\ \hline

  $\alpha_{MIN}$ & the minimum step of the local descent $[10^{-15}, 10^{-1}]$, equal to $10^{-12}$ (Var2).\\

\end{tabularx}
\end{table}

\begin{algorithm}[!ht]
	\caption{Raider method}
  \label{alg:Raider}
	\begin{algorithmic}[1]		
  \REQUIRE algorithmic parameters setup, initial $x^{0}$.
		
		\FOR{$K = 0, \hdots, T_{\text{out}}$}
		\IF	{$ || \nabla f(x^{K}) || \le \varepsilon _{NG}$}
		  \STATE $x^* =  x^{K}$
		  \STATE \textbf{return}
		\ENDIF		
		
		\STATE $M_{GRAD} = \max \{ | \nabla_i f(x^{K}) |,\; i=\overline{1,n}\} $.
		
		\FOR{$i = 1, \hdots, n$}
		\IF {$| \nabla_i f(x^{K}) | >D_{LEVEL} \cdot M_{GRAD} $}
		\STATE  $\mathbf{y}_{i}^{K} =x^{K} _{i} -\nabla_i f(x^{K} )$
		\ELSE 
		\STATE $\mathbf{y}_{i}^{K} =x^{K} _{i} $.
		\ENDIF
		\ENDFOR
		
		\STATE  $\alpha =1$. 
		\STATE $x(\alpha )=x^{K} +\alpha (\mathbf{y}^{K} -x^{K})$ \label{raider:st2}
		
		\IF {$f(x(\alpha))<f(x^{K} )$} 
		  \STATE $x^{K+1} =x(\alpha )$
		\ELSE
		  \STATE  $\alpha =\alpha / K_{\alpha} $.
  		\IF {$\alpha \le \alpha _{MIN}$}
       \STATE $x^* =  x^{K}$
       \STATE \textbf{return}
		  \ELSE
		    \STATE \textbf{go to line~\ref{raider:st2}}
		  \ENDIF
    \ENDIF

	\ENDFOR
	\ENSURE $x^*$ if the early termination condition was reached, otherwise the final iterate $x_{T_{\text{out}} + 1}$.
	\end{algorithmic}
\end{algorithm}

\subsection{Levenberg--Marquardt Method}

It was also explored other ways to select the $S(x^k)$ matrix, whose scheme looks like $S(x^k) = \lambda^k B(x^k)$, where $B(x^k)$ selected as described in table~\ref{tbl:lm_variants}.

\begin{table}[H]
\centering
\caption{Levenberg--Marquardt method variants}
\label{tbl:lm_variants}
\begin{tabular}{r|l}
Variant & $B(x^k)$ \\ \hline
  1 & $I = \diag\{1,1,...,1\}$ \\
  2 & $\diag\{G\} = \diag\{G_{11}(x^k), G_{22}(x^k),...,G_{nn}(x^k)\}$, where $G(x^k) = J(x^k)^T J(x^k)$ \\
  3 & $\diag\{\sqrt{G}\} = \diag\{\sqrt{G_{11}(x^k)}, \sqrt{G_{22}(x^k)},...,\sqrt{G_{nn}(x^k)}\}$
\end{tabular}
\end{table}


Algorithm~\ref{alg:LM} describes a modification of the Levenberg--Marquardt method implemented for the presented experiments. Parameter values used in computational experiments: $\alpha = 0.5$, $\beta = 10$, $\lambda_0 = 10^{-1}$.

\begin{algorithm}
\caption{Levenberg--Marquardt method}
\label{alg:LM}
\begin{algorithmic}[1]
  \REQUIRE $0 < \alpha < 1$, $\beta > 1$, $\lambda^0 > 0$
  \FOR {$k = 1,\dots,N$}
    \STATE $\lambda^k = \lambda^{k-1}$
    \STATE Calc $J(x^k)$, $B(x^k)$
    \STATE Try to find direction $d^k$ (solve linear system with Cholesky decomposition for example):\label{lm:solve}
    \[
      \left( J(x^k)^T J(x^k) + \lambda^k B^k \right) d^k = -J(x^k)^T f(x^k)
    \]
    \IF {$d^k$ is not found (linear system cannot be solved)}
      \STATE $\lambda^k = \lambda^k \cdot \beta$, \textbf{goto line~\ref{lm:solve}}
    \ENDIF
    \STATE $x^{k+1} = x^k + d$
    \IF {$F(x^{k+1}) \ge F(x^k)$}
      \STATE $\lambda^k = \lambda^k \cdot \beta$, \textbf{goto line~\ref{lm:solve}}
    \ELSE
      \STATE $\lambda^k = \lambda^k \cdot \alpha$
    \ENDIF
  \ENDFOR
\end{algorithmic}
\end{algorithm}

\subsection{Modification of B.T. Polyak Method}

\begin{table}[H]
\centering
\caption{B.T. Polyak Method algorithmic parameters}
\begin{tabularx}{\linewidth}{l|X}
  $f^*$ & the lower bound for the optimal value of the function from $[-10^{10}, 10^{10}]$, equal to $0.0$. \\ \hline

  $\varepsilon_{NG}$ & the accuracy of the stopping criterion according to the gradient norm from $[10^{-12}, 10^2]$, equal to $10^{-5}$.\\ \hline



  $\alpha_{MIN}$ & the minimum value of local step $[10^{-15}, 10^{-1}]$, $10^{-12}$.\\ \hline

  $K_{\alpha}$ & the coefficient of the step compression $[1.1, 10.0]$, $3.0$.\\ \hline
\end{tabularx}
\end{table}

\begin{algorithm}[!ht]
	\caption{Polyak method}
  \label{alg:PLK}
	\begin{algorithmic}[1]
		\REQUIRE initial $\mathbf{x}^{0}$, number of iterations $T_{\text{out}}$.
		\FOR{$K = 0, \hdots, T_{\text{out}}$}\label{st1}
		\IF	{$\| \nabla f(x^{K} ) \| \le \varepsilon_{NG}$}
		  \STATE $x^* = x^{K}$
  		\STATE \textbf{return}
		\ENDIF
		
		\IF{Variant1}
  		\STATE $\alpha = (f(x^{K}) - f_{NIZ}) / || \nabla f(x^K)||^2$
		\ENDIF
  		\IF{Variant2}
      \STATE $\alpha = 2 (f(x^{K}) - f_{NIZ}) / || \nabla f(x^K)||^2$
		\ENDIF

		\STATE $x^{K+1} = x^K - \alpha \nabla f(x^K)$\label{plk:st2}
 
		\IF {$f(x^{K+1}) \ge f(x^K) $}
      \STATE $\alpha = \alpha / K_{\alpha}$
      \IF {$\alpha < \alpha_{MIN}$}
  		  \STATE $x^* = x^{K}$
    		\STATE \textbf{return}
      \ELSE
        \STATE \textbf{go to line \ref{plk:st2}}
      \ENDIF
		\ENDIF
  \ENDFOR
  \ENSURE $\mathbf{x}^*$ if the early termination condition was reached, otherwise the final iterate $x_{T_{\text{out}} + 1}$.
	\end{algorithmic}
\end{algorithm}

\end{document}